\documentclass[12pt]{l4dc2022}

\usepackage{algorithm2e}
\usepackage{comment}
\usepackage{amsmath}
\usepackage{algcompatible}
\usepackage{etoolbox}

\usepackage{algorithm}

\AtBeginEnvironment{algorithm2e}{\algoequations}
\AtEndEnvironment{algorithm2e}{\restoreequations}
\newcounter{algosavedequation}
\newcommand{\algoequations}{%
  \setcounter{algosavedequation}{\value{equation}+1}%
  \setcounter{equation}{0}%
  \renewcommand{\theequation}{\arabic{algosavedequation}\alph{equation}}
}
\newcommand{\restoreequations}{%
  \setcounter{equation}{\value{algosavedequation}}%
}

\allowdisplaybreaks[3]

\newtheorem{assumption}{Assumption}
\newtheorem{problem}{Problem}
\newtheorem{rem}{Remark}

\newcommand{\real}{\mathbb{R}} 

\newcommand{\naturalset}{\mathbb{N}}
\newcommand{\naturalpos}{\mathbb{N}_{>0}}
\newcommand{\realpos}{\mathbb{R}_{> 0}}
\newcommand{\realnneg}{\mathbb{R}_{\geq 0}}

\newcommand{\tsp}{\mathsf{T}} 
\newcommand{\pinv}{\dagger} 
\newcommand{\inv}{{\negat 1}} 
\newcommand{\negat}{\scalebox{0.75}[.9]{\( - \)}}

\newcommand*{\QEDB}{\hfill\ensuremath{\square}}
\newcommand*{\QEDBL}{\hfill\ensuremath{\blacksquare}}
\newcommand\oprocendsymbol{\hbox{$\square$}}
\newcommand\oprocend{\relax\ifmmode\else\unskip\hfill%
\fi\oprocendsymbol}

\newcommand{\R}{\mathbb{R}} 

\newcommand{\map}[3]{#1: #2 \rightarrow #3}

\newcommand{\sbs}[2]{{#1}_{\textup{#2}}}

\newcommand{\until}[1]{\{1,\dots,#1\}}

\newcommand{\ess}{\operatorname{ess}}

\newcommand{\norm}[1]{\Vert #1 \Vert}

\title[Data-Enabled  Gradient Flow as Feedback Controller]{Data-Enabled  Gradient Flow as Feedback Controller: Regulation of  Linear Dynamical Systems to Minimizers of Unknown Functions}

\usepackage{times}

 \author{\Name{Liliaokeawawa Cothren} \Email{liliaokeawawa.cothren@colorado.edu}\\
  \Name{Gianluca Bianchin} \Email{gianluca.bianchin@colorado.edu}\\
  \Name{Emiliano Dall'Anese} \Email{emiliano.dallanese@colorado.edu}\\
  \addr Department of Electrical, Computer, and Energy Engineering, University of Colorado Boulder}

\begin{document}

\maketitle

\begin{abstract}%

This paper considers the problem of regulating a linear dynamical system to the solution of a convex optimization problem with an unknown or partially-known cost. We design a data-driven feedback controller -- based on gradient flow dynamics -- that (i) is augmented with learning methods to estimate the cost function based on infrequent (and possibly noisy) functional evaluations; and, concurrently, (ii) is designed to drive the inputs and outputs of the dynamical system to the optimizer of  the problem. We derive sufficient conditions on the learning error and the controller  gain to ensure that the error between the optimizer of the problem and the state of the closed-loop system is ultimately bounded; the error bound accounts  for the functional estimation errors and the temporal variability of the unknown  disturbance affecting the linear dynamical system. Our results directly lead to  exponential input-to-state stability of the closed-loop system. The proposed method and the theoretical bounds are validated numerically.

\end{abstract}

\begin{keywords}%
 Learning-based control; learning-based optimization; gradient flow; output regulation. 
\end{keywords}

\section{Introduction}
\label{S:introduction}

In this paper, we consider the problem of designing feedback controllers to steer the output of a  linear time-invariant (LTI) dynamical system towards the solution of a convex optimization problem with unknown costs. The design of  controllers inspired by optimization algorithms has received attention recently; see, e.g.,~\cite{Jokic2009controller, brunner2012feedback,lawrence2018optimal,hauswirth2020timescale,MC-ED-AB:20,zheng2019implicit,bianchin2020online} and the recent survey by~\cite{hauswirth2021optimization}. These methods have been utilized to solve control problems in, e.g., power systems in~\cite{Hirata,menta2018stability}, transportation systems in~\cite{bianchin2021time}, robotics in~\cite{zheng2019implicit}, and epidemics in~\cite{bianchin2021planning}. 

A common denominator in the works mentioned above  is that the cost of the optimization problem associated with the dynamical system is \emph{known}, and first- and second-order information is easily accessible at any time. One open  research  question is  whether controllers  can  be  synthesized when the cost of the  optimization problem is  \emph{unknown} or \emph{partially known}. Towards this direction, in this paper we consider unconstrained convex optimization problems with unknown costs associated with the LTI dynamical systems. We investigate the design  of  \emph{data-driven} feedback controllers based on online gradient flow dynamics that: \emph{(i)} leverage learning methods to estimate the cost function based on infrequent (and possibly noisy) functional evaluations; and \emph{(ii)} are designed to concurrently drive the inputs and outputs of the dynamical system to the optimizer of  the problem within a bounded error. Our learning procedure hinges on a basis expansion for the cost and leverages methods such as least-squares, ridge regression, and sparse linear regression; see~\cite{hastie2009elements,tibshirani1996regression}. In addition, our results are also directly applicable to  cases where residual neural networks are utilized to estimate the cost; see \cite{tabuada2020universal}.

We consider a cost that includes the sum of a loss associated with the inputs and a loss function associated with the outputs. Both costs may be unknown or  parametrized by unknown parameters; we assume that functional evaluations are provided at irregular intervals, due to underlying communication or processing bottlenecks (see, for example, delays in power system metering systems in~\cite{luan2013data} and in transportation systems in~\cite{gundling2020efficient}). As an example of a function with unknown parameters associated with the outputs, take $\psi(y) = \|y - r\|^2$, where $\map{y}{\realpos}{ \mathbb{R}^p}$ is the system output and $r \in \mathbb{R}^p$ is unknown to the controller.  Additional examples include cases where $\psi(\cdot)$ represents a barrier function associated with unknown sets; see, e.g.,~\cite{robey2020learning,taylor2020learning}. Regarding the function associated with the inputs,  another scenario where it is unknown is when it  captures objectives of  users interacting with the system; see~\cite{simonetto2021personalized, notarnicola2021distributed,fabiani2021learning,ospina2020personalized,luo2020socially}. In this case, the  loss  models objectives such as dissatisfaction, discomfort, etc. For example, in a platooning problem, the control input is represented by the speed of the vehicles, and the loss function captures  the sense of safety for drivers~\cite{Nunen}. In lieu of synthetic models (that may not represent accurately the user's objectives~\cite{munir2013cyber,bourgin2019cognitive}), one learns the loss based on evaluations infrequently provided by the user.

\textbf{Related Works.} We note that a key differentiating aspect relative to extremum seeking methods (see, e.g.,~\cite{krstic2000stability,ariyur2003real,teel2001solving} and many others), the  Q-learning of~\cite{devraj2017zap}, and methods based on concurrent learning~\cite{chowdhary2010concurrent, chowdhary2013concurrent,poveda2021data} is that we consider  a setting where only sporadic functional evaluations are  available (i.e., we do not have continuous access to functional evaluations). Regarding the problem of regulating  LTI systems  towards solutions of optimization problems, existing approaches leveraged gradient flows in~\cite{menta2018stability,bianchin2020online}, proximal-methods in~\cite{MC-ED-AB:20}, saddle-flows in~\cite{brunner2012feedback}, prediction-correction methods in~\cite{zheng2019implicit}, and the hybrid accelerated methods proposed in~\cite{bianchin2020online}. Plants with (smooth) nonlinear dynamics were considered  in~\cite{brunner2012feedback,hauswirth2020timescale}, and switched LTI systems in~\cite{bianchin2021time}. A joint stabilization and regulation problem was considered in~\cite{lawrence2018linear,lawrence2018optimal}. See also the recent survey by~\cite{hauswirth2021optimization}. In all these works, the cost function is assumed to be known; here, we tackle the problem of jointly learning the cost and performing the regulation task.    

Our setup is aligned with~\cite{simonetto2021personalized, ospina2020personalized}, where  Gaussian Processes are utilized to learn cost functions based on infrequent functional evaluations, and \cite{notarnicola2021distributed}, where the cost is estimated via recursive least squares method. However, these works focus on discrete-time algorithms and, more importantly, have no dynamical system implemented in closed-loop with the algorithms. 
Few recent works considered controllers that are learned using neural networks; see, e.g.,~\cite{karg2020stability,HY-PS-MA:21,Marchi2022}, and the work on reinforcement learning in~\cite{MJ-JL:20}. With respect to this literature, we utilize learning methods to estimate the cost, and we use a gradient-flow controller based on the estimated cost.  
Finally, similarly to~\cite{sontag2021remarks}, we study the input-to-state 
stability (ISS) property of perturbed gradient flows. Differently 
from~\cite{sontag2021remarks}, in this work we consider the interconnection 
between a perturbed gradient-flow and an LTI system.

\textbf{Contributions.} Our contribution is threefold.  \emph{(C1)} We design a data-driven feedback controller to steer the inputs and outputs of an LTI system towards the optimizer of a convex optimization problem; the controller does not require knowledge of the unknown and time-varying exogenous inputs affecting the system. The controller leverages  methods that  learn the cost functions of the optimization problem from historical information (as a starting estimate) and through infrequent functional evaluations during the operation of the controller. Our setting accounts for cases where the cost function admits a representation through a finite set of basis functions, and the more general case where we approximate the function using a truncated basis expansion. 
\emph{(C2)} We leverage singular-perturbation arguments (as in~\cite[Ch.~11]{Khalil:1173048} and in, e.g.,~\cite{hauswirth2020timescale,bianchin2020online}) and the theory of perturbed systems~\cite[Ch.~9]{Khalil:1173048} to derive sufficient conditions
on the learning error and the controller  gain to ensure that the error 
between the optimizer of the problem and the state of the closed-loop system 
is ultimately bounded. \emph{(C3)}   We verify the stability claims and the analytical bounds through a representative set of simulations. 

\textbf{Organization.} Section \ref{S:problemFormulation} outlines the problem formulation and the main assumptions; Section~\ref{S:mainResults} presents the main data-driven control framework and the stability results. Representative numerical simulations are
presented in Section~\ref{S:numericalVerification}, and Section~\ref{S:conclusion} concludes the paper. The proofs of the main results are reported in the Appendix\footnote{\textbf{Notation.}
We denote by \( \mathbb{N}, \naturalpos, \R, \realpos, \text{ and } \realnneg \) the set of natural numbers, the set of positive natural numbers, the set of real numbers, the set of positive real numbers, and the set of non-negative real numbers, respectively. For vectors \( x \in \R^n \) and \( u \in \R^m \),
\( \| x \| \) denotes the  Euclidean norm of $x$ and \( (x,u) \in \R^{n + m}\) denotes their vector concatenation. If  \( x,u \in \R^n \), then \( (x^\top ,u^\top ) = [x^\top; u^\top] \in 
\R^{2 \times n}\) denotes the matrix with rows given by 
$x^\top$ and $u^\top$. 
For a symmetric matrix $W \in \real^{n \times n}$, 
$W \succ 0$ denotes that $W$ is positive definite and $W \succeq 0$ denotes 
that $W$ is positive semidefinite. Moreover, we let 
\( \bar{\lambda}(W) \) and \( \underline{\lambda}(W) \) denote 
the largest and smallest eigenvalues of $W$, respectively. 
For a continuously differentiable function 
$\map{\phi}{\real^n}{\real}$, we denote its gradient by 
\( \nabla \phi(x) \in \real^n\). }.

\section{Problem Formulation}\label{S:problemFormulation}

We consider continuous-time linear  dynamical systems described by:
\begin{align}
\label{eq:plantModel}
\dot x &=  A x + B u + E w_t, & 
y &= C x + D w_t,
\end{align}
where $\map{x}{\realnneg}{\real^n}$ is the state, 
$\map{u}{\realnneg}{\real^m}$ is the input, 
$\map{y}{\realnneg}{\real^p}$ is the output,
 $\map{w_t}{\realnneg}{\real^q}$ is an unknown 
and time-varying exogenous input or disturbance, and $A, B, C, D$, and
$E$ are matrices of appropriate dimensions. 
We make the following assumptions on~\eqref{eq:plantModel}.
\begin{assumption}
\label{ass:stabilityPlant}
The matrix $A$ is Hurwitz stable; namely, for any 
$Q \in \real^{n\times n}, Q \succ 0$, there exists 
$P \in \real^{n \times n}, P \succ 0,$ such that 
$A^\top P +PA = -Q$.
\QEDB
\end{assumption}
\vspace{-.4cm}
\begin{assumption}
\label{as:disturbance}
The function \( t \mapsto w_t\) is locally absolutely continuous. 
\QEDB
\end{assumption}
\vspace{-.2cm}
Under Assumption \ref{ass:stabilityPlant}, for given vectors 
$\sbs{u}{eq} \in \real^m$ and $\sbs{w}{eq} \in \real^q$, 
\eqref{eq:plantModel} has a unique stable equilibrium point 
$\sbs{x}{eq} =-A^{-1}(B \sbs{u}{eq}+E \sbs{w}{eq})$; see, e.g.,~\cite[Theorems 4.5 and 4.6]{Khalil:1173048}. 
Moreover, at equilibrium, the relationship between system 
inputs and outputs is given by the algebraic map 
$\sbs{y}{eq} = G \sbs{u}{eq} + H \sbs{w}{eq}$, where 
$G := -C A^\inv B$ and $H := D -C A^\inv E$. 
Assumption~\ref{as:disturbance} characterizes how the exogenous inputs can 
vary over time. 
 
We consider the problem of developing a feedback controller, inspired by 
online optimization methods, to regulate \eqref{eq:plantModel} to the 
solutions of the following time-dependent optimization problem:
%
\begin{align}
\label{opt:objectiveproblem}
u_t^* \in  
\underset{\bar u \in \real^m}{\arg \min}  
~~~ &  \phi (\bar u) + \psi (G  \bar u  +H w_t),
\end{align}
for all $t \in \realnneg$, where $\map{\phi}{\real^m}{\real}$ and $\map{\psi}{\real^p}{\real}$ are cost functions associated with the system's inputs and outputs, respectively. The optimization problem \eqref{opt:objectiveproblem} formalizes an equilibrium selection problem  for which the objective is to select an optimal 
input $u^*_t$ for the system \eqref{eq:plantModel} (and, 
consequently, the corresponding steady-state output 
$y^*_t = G u^*_t + H w_t$) that minimizes the cost specified by
the the loss functions $\phi(\cdot)$ and $\psi(\cdot)$. 
We note that, since the cost function is parametrized by $w_t$, the solutions 
of~\eqref{opt:objectiveproblem} are also time-varying, and thus
define optimal trajectories (the sub-script $t$ is utilized to emphasize the temporal variability of $w_t$ and, consequently, that of $u^*_t$). 

\begin{rem}{\bf \textit{(Relationship with output regulation)}}
The problem \eqref{opt:objectiveproblem} formalizes an 
optimal regulation problem with steady-state constraints similar 
to the well-established output-regulation problem~\cite{ED:76}; 
with respect to the classical framework, in our setting the optimal
trajectories are not generated by an exosystems (i.e., a known 
autonomous linear model) but instead are specified as the solution of a convex optimization problem.
\QEDB
\end{rem}

In this work, we focus on a setting where the exogenous input $w_t$ is \emph{unknown} and the cost functions $u \mapsto \phi(u)$ and $y \mapsto \psi(y)$ are  \emph{unknown} as explained  in Section~\ref{S:introduction}. 
In this setup, the output regulation problem tackled in this paper is summarized as follows. 

\begin{problem} 
Design a data-driven output-feedback controller 
for~\eqref{eq:plantModel}  that learns the cost functions $u \mapsto \phi(u)$ and $y \mapsto \psi(y)$ from infrequent functional evaluations while concurrently driving the inputs and outputs of \eqref{eq:plantModel} to the time-varying optimizer of 
\eqref{opt:objectiveproblem} up to an error that accounts 
for the functional estimation errors and the temporal variability of the unknown 
disturbance. 
\end{problem}


Although unknown, we impose the following regularity 
assumptions on the cost functions.
\begin{assumption}
\label{as:lipschitzPHI}
The function $u \mapsto \phi(u)$ is 
continuously-differentiable, convex, and $\ell_u$-smooth, for 
some $\ell_u \geq 0$; namely, $\exists~\ell_u \geq 0$ such that
$\|\nabla \phi(u) - \nabla \phi(u')\| \leq \ell_u \|u-u'\|$
holds $\forall \,\, u, u' \in \R^m$.
\QEDB
\end{assumption}
\vspace{-.2cm}
\begin{assumption}
\label{as:lipschitzPSI}
The function $y \mapsto \psi(y)$ is 
continuously-differentiable, convex, and $\ell_y$-smooth, for 
some $\ell_y \geq 0$; namely, $\exists~\ell_y \geq 0$ such that
$\|\nabla \psi(y) - \nabla \psi(y')\| \leq \ell_y \|y-y'\|$
holds $\forall \,\, y, y' \in \R^p$.
\QEDB
\end{assumption}
\vspace{-.3cm}
\begin{assumption}
\label{as:Strconvex}
For any $w \in \real^q$, the composite cost
$u \mapsto \phi(u) + \psi(Gu + H w)$ is  $\mu_u$-strongly 
convex, with $\mu_u > 0$.
\QEDB
\end{assumption}
\vspace{-.1cm}

Assumptions~\ref{as:lipschitzPHI}-\ref{as:lipschitzPSI} imply 
that $u \mapsto \phi(u) + \psi(Gu + H w_t)$ is $\ell$-smooth, 
with $\ell:= \ell_u+\|G\|^2 \ell_y \geq 0$. 
Two implications follow from Assumption~\ref{as:Strconvex}: (i) the optimizer 
$u_t^*$ is unique, and (ii)  $\phi(u) + \psi(Gu + H w_t)$ satisfies the 
Polyak-\L ojasiewicz (PL) inequality as shown in~\cite{karimi2016linear}; 
namely, 
\begin{align*}
\Vert \nabla \phi (u) + G^\top \nabla \psi(G u + H w_t) \|^2 
\geq 
2 \mu_u (\phi(u) + \psi(G u + H w_t) - \phi(u^*_t) - \psi(G u_t^* + H w_t)),
\end{align*}
holds for all $u \in \R^m$. Regarding the functions $\phi(u)$  and $\psi(y)$, 
we make the following assumptions.
\begin{assumption}
\label{as:Phi}
The function \(u \mapsto \phi(u)\) admits the 
representation \(\phi(u) =  \sum_{i=1}^{N_b}   \alpha_{i} b_i(u)\),
for some $N_b \in \naturalpos \cup \{+\infty\}$, where for all $i \in \until {N_b}$, 
\(b_i: \mathbb{R}^m \rightarrow \mathbb{R}\) are continuously  
differentiable  basis functions and  \(\alpha_{i} \in \R\) are fitting parameters.
\QEDB
\end{assumption}
\vspace{-.2cm}
\begin{assumption}
\label{as:Psi}
The function \(y \mapsto \psi(y)\) admits the  
representation \(\psi(y) =  \sum_{i=1}^{M_b}  \rho_{i} d_i(y)\),
for some $M_b \in \naturalpos \cup \{+\infty\}$, where for all $i \in \until {M_b}$,  
$d_i: \mathbb{R}^p \rightarrow \mathbb{R}$ are continuously  
differentiable  basis functions and \(\rho_{i} \in \R\) are 
fitting parameters.
\QEDB
\end{assumption}
\vspace{-.2cm}

In compact form, we denote $\alpha := (\alpha_1, \dots , \alpha_N)$,
$\rho := (\rho_1, \dots , \rho_M)$, $b(u):= (b_1(u), \ldots, b_N(u))$, and 
$d(u):= (d_1(u), \ldots, d_M(u))$. Moreover, we let $\nabla b$ denote the 
Jacobian of $b(u)$, and $\nabla d$  the Jacobian of $d(u)$. 
We illustrate  the above  assumptions through the following two examples.

\vspace{-.2cm}

\begin{example} 
\label{ex:nonparameteric}
\emph{\bf \textit{(Non-parametric models)}}. A representation as in Assumptions~\ref{as:Phi}-\ref{as:Psi} can be obtained by using tools from Reproducing Kernel Hilbert Spaces, in which a function $\phi$ defined over a measurable space is estimated via interpolation based on symmetric, positive definite kernel functions~\cite[Chapter 5]{hastie2009elements},~\cite{bazerque2013nonparametric}. Additional non-parametric models utilize orthonormal basis functions such as polynomials, or can leverage radial basis functions and multilayer feed-forward networks~\cite{hornik1989multilayer}.
\vspace{-.2cm}
\end{example}

\begin{example} 
\label{ex:parameteric}
\emph{\bf \textit{(Convex parametric models)}}. Consider
the  cost  
$\phi(u) = \frac{1}{2} u^\top \Upsilon u + \upsilon^\top u + r$, where 
$\Upsilon \in \R^{m \times m}$, $\Upsilon \succeq 0$, $\upsilon \in \mathbb{R}^m$, 
$r \in \R$. Taking as an example $m=2$, $\Upsilon=[\Upsilon_{ij}]$, $\upsilon=(c_1,c_2)$ 
and $ u = (u_1, u_2)$, the function admits the representation 
 in Assumption \ref{as:Phi} with $N_b=6$, 
where 
$b_1(u) = 1$, 
$b_2(u) = u_1$, 
$b_3(u) = u_2$, 
$b_4(u) = u_1^2/2$, 
$b_5(u) = u_1 u_2$, and 
$b_6(u) = u_2^2/2$, and 
$\alpha =(r, \upsilon_1, \upsilon_2, \Upsilon_{11}, \Upsilon_{12}, \Upsilon_{22})$ (with the constraint $\Upsilon_{12} = \Upsilon_{21}$). 

\vspace{-.2cm}
\end{example}

In our learning strategy, we consider a finite number of 
basis functions as $\tilde \phi(u) = \sum_{i=1}^N  \alpha_i b_i(u)$ and 
$\tilde \psi(y) = \sum_{j=1}^M \rho_j d_j(y)$, where $N < \infty$ and 
$M < \infty$; in particular, we assume that $N\leq N_b$ and $M \leq M_b$. Accordingly, we
will consider two cases:

\noindent \emph{\textbf{(Case~1)}} The functions $\phi$ and $\psi$ are  
represented by a finite number of basis functions (i.e., $N_b < \infty$ and 
$M_b < \infty$), and we utilize $N = N_b$ and $M = M_b$.

\noindent \emph{\textbf{(Case~2)}} At least one of the functions $\phi$ or 
$\psi$ is represented by an infinite number of basis functions (i.e., 
$N_b = \infty$ and/or $M_b = \infty$). In this case, we utilize $N <\infty$ 
and $M < \infty$.

We remark that the scenario in~\textit{(Case 2)} can be of interest also in 
cases where $N_b$ or $M_b$ are finite but sufficiently large making the 
learning technique computationally intractable. In this
case, in our framework we will approximate the two functions by utilizing 
$N \ll N_b$ and/or $M \ll M_b$.
Finally, we note that an instance of \emph{(Case 1)} is the convex parametric 
model in Example~\ref{ex:parameteric}. For the latter, for a general nonparametric model as in Example~\ref{ex:nonparameteric}, due  to the Weierstrass high-order approximation theorem, the error in the representation of the function over a compact set goes to zero  with the increasing of  $N_b$ and $M_b$; however, a limited number of basis functions may be selected for model complexity considerations (see, e.g.,~\cite{hastie2009elements} and ~\cite{bazerque2013nonparametric}).

\section{Gradient-flow Controller with Concurrent Learning}\label{S:mainResults}

To address Problem~1, we propose a concurrent learning and online 
optimization scheme in which, at every time $t$, we: 
(i)~process a newly available functional evaluation $(u, \phi(u))$ (resp. $(y,\psi(y))$) via a learning method to 
determine an estimate $\hat \alpha_t$ of 
the vector of parameters $\alpha$ (resp. to determine an estimate 
$\hat \rho_t$ of the vector $\rho)$ at time $t$;  and
(ii) use an approximate gradient-flow  parametrized by the current estimates $(\hat \alpha_t, \hat \rho_t)$ to steer~\eqref{eq:plantModel} towards the optimal solution of the 
 problem~\eqref{opt:objectiveproblem}.

For the learning procedure, we utilize both noisy functional evaluations 
received during the operation of the algorithm as well as from historical 
data.  
Accordingly, let $\hat{\phi}_{t_k} = \phi(u_{t_k}) + \varepsilon_{t_k}$ be the
noisy functional evaluation received at time $t_k \in \R_{\geq 0}$, 
$k \in \mathbb{N}$, at the point $u_{t_k} \in \mathbb{R}^m$, in the presence 
of measurement noise $\varepsilon_{t_k} \in \mathbb{R}^m$. 
(By convention, we let $t_k = 0$, $k \in \mathbb{N}$, if the $k$-th pair 
$\{(\hat{\phi}_{t_k}, u_{t_k})\}$ is derived from historical data.) 
Similarly, we let $\hat{\psi}_{t_j} = \psi(y_{t_j}) + v_{t_j}$ be a noisy 
functional evaluation of $\psi(y)$, received at time $t_j \in \R_{\geq 0}$, 
$j \in \mathbb{N}$, in the presence of output noise 
$v_{t_j} \in \mathbb{R}^p$. We consider learning methods that yield estimates of $\alpha$ and $\rho$ at time $t$ based on the data  
$\{\hat{\phi}_{t_k}, u_{t_k}\}_{t_k \leq t}$ and 
$\{\hat{\psi}_{t_j}, y_{t_j}\}_{t_j \leq t}$, respectively.   

\begin{figure}[t!]
\centering 
\includegraphics[width=.9\columnwidth]{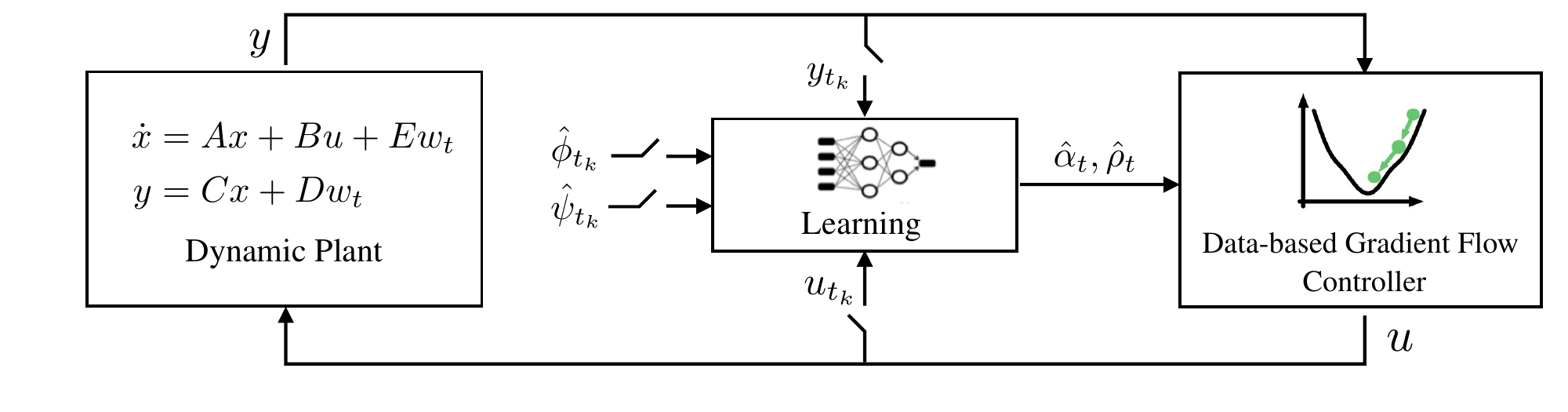}
\vspace{-.4cm}
\caption{Block diagram of the proposed gradient flow controller with concurrent learning.}
\label{fig:system}
\vspace{-.2cm}
\end{figure}

\begin{algorithm2e}[h!]\label{eq:plantModelandController}
\caption{Data-based online gradient-flow controller}
\DontPrintSemicolon
\KwIn{
$\hat \alpha_{t_0}, \hat \rho_{t_0}$ based on recorded data;  gain $\eta > 0$, number of basis functions $N$, $M$. \\
}
\vspace{0.15cm}
\textbf{For all} $t \geq t_0$: 

\textbf{\# Learning }\\
\uIf{ new data point $(\hat \phi(u_t),u_t)$ obtained 
}{
\vspace{-0.5cm}
\begin{flalign}\label{eq:plantModelandController-alpha}
    \hat \alpha_{t} &\gets \operatorname{parameter-learning}(\{(\hat{\phi}_{t_k}, u_{t_k})\}_{t_k \leq t}) &
\end{flalign}
}
\If{ new data point $(\hat \psi(y_t),y_t)$ obtained}{
\vspace{-0.5cm}
\begin{flalign}\label{eq:plantModelandController-rho}
    \hat \rho_{t} &\gets \operatorname{parameter-learning}(\{(\hat{\psi}_{t_j}, y_{t_j})\}_{t_j \leq t}) &
    \vspace{-.3cm}
\end{flalign}
}
\textbf{\# Feedback Control }\\
\vspace{-0.75cm}
\begin{flalign}\label{eq:plantModelandController-b}
    \dot x &= A x + B u + E w_t , \,\, y  = C x + D w_t \nonumber &\\
    \dot u &= -\eta ( \nabla b(u)^\top \hat \alpha_t + G^\top \nabla d(y)^\top \hat{\rho}_t). &
\end{flalign}
\vspace{-0.5cm}
\vspace{-.2cm}
\end{algorithm2e}
The proposed scheme is illustrated in Figure~\ref{fig:system}, and it is described by the pseudo-code in Framework~\ref{eq:plantModelandController} 
(with $t_0$ denoting the initial time).
According to the proposed method, an initial estimate of the parameters $\alpha$ and $\rho$ is obtained using recorded data, where   $\operatorname{parameter-learning}(\cdot)$ is a  map that represents an estimation step;  examples of learning methods will be provided in Section~\ref{sec:estimation}. Since for any time $t \geq t_0$ the parameters are updated only when a new functional evaluation is available, we let $t \mapsto \hat{\alpha}_t$ and $t \mapsto \hat{\rho}_t$ be piece-wise constant right-continuous functions that represent the most up-to-date estimates of $\alpha$ and $\rho$, respectively. We utilize a gradient-flow controller, as shown in~\eqref{eq:plantModelandController-b}, where the gain $\eta >0$ induces a time-scale separation between the plant and the controller. Finally, we notice that since the vector field characterizing 
\eqref{eq:plantModelandController-b} is piece-wise continuous in time, 
the initial value problem \eqref{eq:plantModelandController-b} always 
admits a local solution that is unique  \cite[Thm 3.1]{Khalil:1173048}.

\subsection{Main results}
\label{sec:mainResults}

In this section, we characterize the transient performance 
of Framework~\ref{eq:plantModelandController}. To this end, let 
\begin{align}
\label{eq:trackingErrorDef}
z &:= (u - u^*_t, x - x^*_t), & 
x^*_t &= - A^\inv B u^*_t - A^\inv E w_t,
\end{align}
denote the error between the state of 
\eqref{eq:plantModelandController-b} and the equilibrium of 
\eqref{opt:objectiveproblem}, respectively. In the following, we provide sufficient conditions on $\eta$ and on the estimation errors $\alpha - \hat{\alpha}_t$, $\rho - \hat{\rho}_t$ so that the  interconnection between the plant and the data-based  controller~\eqref{eq:plantModelandController-b} is exponentially stable. 

The first result is stated for  \emph{\textbf{(Case~1)}}, where the functions $\phi$ and $\psi$ are represented by a finite number of basis functions, and we set precisely $N = N_b$ and $M = M_b$.   

\begin{theorem}{\bf \textit{(Control bound for finite number of basis functions)}}
\label{thm:finiteBasisFunctions} 
Let Assumptions \ref{ass:stabilityPlant}-\ref{as:Psi}
be satisfied, and assume that $N = N_b < \infty$ and $M = M_b < \infty$. Let $z(t)$ be defined as in \eqref{eq:trackingErrorDef}, with $(u(t), x(t))$ the state of 
\eqref{eq:plantModelandController-b}.
Suppose that the learning errors satisfy,
\begin{align}
\label{eq:epsilonUpperBound}
 0 < \epsilon := \ell_u \sup_{t_0 \leq \tau \leq t}\|\alpha - \hat \alpha_\tau\| + \ell_y \|G\|^2 \sup_{t_0 \leq \tau \leq t}\|\rho - \hat \rho_\tau\| < \frac{c_0}{c_3},
\end{align}
where
$c_{0} := s \min \left\{ 2 \mu_u \eta, \underline\lambda (Q) /\bar \lambda (P) \right\}, 
c_{3} :=  \eta \max \left\{2 \ell \mu_u^{-1}, 4 c_1^{-1} \|P A^{-1}B \|\right\},$ $c_1 := \min \{(1 - \theta)\mu_{u}/2,\theta \underline{\lambda}\left(P\right) \},
\theta := \ell_y \|G\|\|C\|/(\ell_y \|G\|\|C\| + 2\|PA^{-1}B\|),$ and $s \in (0,1)$. Suppose further that the controller gain satisfies, 
\begin{equation}
\label{eq:conditiongain}
    0 < \eta < \frac{(1 - s)^2 \underline{\lambda}(Q)}{\left( 2 - s \right)  2\|PA^{-1}B\|\ell_y \|G\|\|C\|}.
\end{equation} 
Then there exists $\kappa_1, \kappa_2,\kappa_3 > 0$ such 
that the error $z(t)$ satisfies
\begin{equation}\label{eq:thmResult}
   \| z(t) \| \leq \kappa_1 e^{- \frac{1}{2} a (t - t_0)}\|z (t_0)\| 
+ \kappa_2 \int_{t_0}^t e^{- \frac{1}{2} a (t - \tau)} 
   \Delta(\tau) d \tau
+ \kappa_3 \int_{t_0}^t e^{- \frac{1}{2} a (t - \tau)}  
\norm{\dot w_\tau } d \tau,
\end{equation}
for all  $t \geq t_0 \geq 0$, where 
$a := c_0 - \epsilon \, c_3$ and 
$\Delta(\tau) :=  \| \nabla b(u^*_\tau) \| \|\alpha - \hat \alpha_\tau\| + \| \nabla d(y^*_\tau) \| \|\rho - \hat \rho_\tau\| $.
\end{theorem}

Detailed expressions  for the constants $\kappa_1, \kappa_2,\kappa_3$ are provided in Appendix~\ref{S:completeresult}, and the proof of Theorem~\ref{thm:finiteBasisFunctions}  is provided in  Appendix~\ref{S:proofs}. Theorem~\ref{thm:finiteBasisFunctions} asserts that if the worst-case  estimation error for the parameters of the cost functions, captured by $\sup_{t_0 \leq \tau \leq t}\|\alpha - \hat \alpha_\tau\|$ and $\sup_{t_0 \leq \tau \leq t}\|\rho - \hat \rho_\tau\|$, is such that  $\epsilon$ satisfies the condition \eqref{eq:epsilonUpperBound}, then a 
sufficiently-small choice of the controller gain $\eta$ guarantees 
exponential convergence of the state 
of~\eqref{eq:plantModelandController-b} to a neighborhood of the
optimizer $u^*_t$ of~\eqref{opt:objectiveproblem} and the corresponding state $x_t^*$. 
In particular, the error $z(t)$ is ultimately bounded 
by two terms: the first depends on the error $\Delta(t)$, which accounts for the estimation error in the function parameters, and the second depends on 
temporal variability of $w_t$ (which affects  
the dynamics of the plant \eqref{eq:plantModel} making the optimizer of~\eqref{opt:objectiveproblem} time varying). The condition on $\epsilon$ also suggests prerequisites on the ``richness'' of the recorded data in the sense that they must yield a sufficiently small estimation error. 

It is important to notice that 
$a = c_0 - \epsilon \, c_3$, which characterizes the rate of  exponential decay,  
is proportional to smallest value between the rate of the convergence of the 
open-loop plant (found in $c_0$ as the ratio 
$\underline \lambda(Q) /\bar \lambda(P)$ 
\cite[Chapter 4, Theorem 4.10]{Khalil:1173048}) and the strong convexity 
parameter $\mu_u$, characterizing the cost function. 
Moreover, the rate of convergence $a$ is proportional to the 
controller gain $\eta$ (as described by $c_0$ and $c_3$), and
inversely proportional to the worst-case estimation errors of
the parameters of the cost, 
$\sup_{t_0 \leq \tau \leq t}\|\alpha - \hat \alpha_\tau\|$ and $\sup_{t_0 \leq \tau \leq t}\|\rho - \hat \rho_\tau\|$.

\begin{rem}{\bf \textit{(Asymptotic behavior and input-to-state stability)}}
\label{rem:boundedness}
Two important implications follow from the statement of 
Theorem \ref{thm:finiteBasisFunctions} as subcases. 

\noindent $\bullet$ If $\lim_{t \rightarrow \infty} \Delta(t) = 0$ and 
$\lim_{t \rightarrow \infty} \dot w_t = 0$, then \eqref{eq:thmResult} guarantees 
that  $\lim_{t \rightarrow \infty} z(t) = 0$, namely, the state of 
\eqref{eq:plantModelandController-b} converges (exactly) to the optimizer
of \eqref{opt:objectiveproblem}. 
See \cite[Lemma 9.6]{Khalil:1173048}.

\noindent $\bullet$ If 
$\ess \sup_{t_0 \leq \tau \leq t} \norm{\Delta(\tau)} < \infty$ and 
$\ess \sup_{t_0 \leq \tau \leq t} \norm{\dot w_\tau} < \infty$, then 
\eqref{eq:thmResult} guarantees that
\begin{align}
\label{eq:ISS}
\| z(t) \| \leq  \kappa_1 e^{- \frac{1}{2} a (t - t_0)}\|z (t_0)\|  + 2a^\inv ( \kappa_2 \ess \sup_{t_0 \leq \tau \leq t} \norm{\Delta(\tau)} + \kappa_3
\ess \sup_{t_0 \leq \tau \leq t} \norm{\dot w_\tau } ).
\end{align}
It follows that the bound \eqref{eq:ISS} guarantees 
input-to state stability of \eqref{eq:plantModelandController-b} (in the sense
of~\cite{sontag1997output,angeli2003input,sontag2021remarks}) with respect to the inputs $\Delta$ 
and $\dot w_t$.
\QEDB
\end{rem}

We now consider \emph{\textbf{(Case~2)}}, where the estimation procedure is affected by a truncation error; that is, when only 
$N \ll N_b$ and $M \ll M_b$ basis functions are utilized by the controller. Accordingly, consider the approximated functions $\tilde \phi(u) = \sum_{i=1}^N  \alpha_i b_i(u)$ and $\tilde \psi(y) = \sum_{j=1}^M \rho_j d_j(y)$, and define the two truncation errors for $\phi(u)$ and $\psi(y)$ as $e_{\phi}(u) := \sum_{i=N+1}^{N_b} b_i(u) \alpha_i$ and $e_{\psi}(y) :=  \sum_{j=M+1}^{M_b}  d_j(y) \rho_j$, respectively.  We make the following regularity assumptions: 

\begin{assumption}\label{as:tildeLipschitz}
The approximated functions $\tilde \phi(u)$ and $\tilde \psi(y)$ are $\ell_u^N$- and $\ell_y^M$-smooth for constants $\ell_u^N \geq 0$ and $\ell_y^M \geq 0$, respectively. 
\end{assumption}

\begin{assumption}
\label{as:truncationLipschitz}
The truncation terms $e_\phi(u)$ and $e_\psi(y)$ have a Lipschitz-continuous gradient with constants $\ell_u^e \geq 0$ and $\ell_y^e \geq 0$, respectively. 
\end{assumption}

Assumption~\ref{as:tildeLipschitz} is satisfied if the basis functions $\{b_i\}_{i = 1}^N$ and $\{d_i\}_{i = 1}^M$ are strongly smooth. Assumption~\ref{as:truncationLipschitz} is a technical condition that is required in our proof; nevertheless, this condition is satisfied for parametric convex models and several non-parametric models;  see~\cite{hastie2009elements, bazerque2013nonparametric}, and the recent results of~\cite{poveda2021data}.

We next characterize the controller transient performance in the presence of 
truncation errors.

\begin{theorem}{\bf \textit{(Stability with function approximation error)}}
\label{thm:truncation}
Let Assumptions \ref{ass:stabilityPlant}-\ref{as:truncationLipschitz}
be satisfied, and assume that $N < N_b$ and $M < M_b$ basis functions 
are utilized in the learning of $\phi$ and $\psi$, respectively.
Suppose that 
\begin{align}
\label{eq:epsilonUpperBoundInfinite}
 0 < \epsilon^\prime := \ell_u^N \sup_{t_0 \leq \tau \leq t}\|\alpha - \hat \alpha_\tau\| + \ell_y^M \|G\|^2 \sup_{t_0 \leq \tau \leq t}\|\rho - \hat \rho_\tau\| + \ell_u^e + \ell_y^e \|G\|^2 < \frac{c_0}{c_3},
\end{align}
where \(s, c_0, c_1, c_3, \theta \) are given in Theorem \ref{thm:finiteBasisFunctions} and $\eta$ satisfies~\eqref{eq:conditiongain}. Then there exists $\kappa_1, \kappa_2,\kappa_3 > 0$ such 
that the error \eqref{eq:trackingErrorDef} satisfies
\begin{equation}\label{eq:thmResultInfinite}
   \| z(t) \| \leq \kappa_1 e^{- \frac{1}{2} a^\prime (t - t_0)}\|z (t_0)\| 
+ \kappa_2 \int_{t_0}^t e^{- \frac{1}{2} a^\prime (t - \tau)} 
   \Xi(\tau) d \tau
+ \kappa_3 \int_{t_0}^t e^{- \frac{1}{2} a^\prime (t - \tau)}  
\norm{\dot w_\tau } d \tau,
\end{equation}
for all  $t \geq t_0 \geq 0$, where 
$a^\prime = c_0 - \epsilon^\prime \,  c_3 $ and
$\Xi(\tau) :=  \| \nabla b(u^*_t) \| \|\alpha - \hat \alpha_\tau\| + \|G\|\| \nabla d(y^*_t) \| \|\rho - \hat \rho_\tau\| + \|\nabla e_{\phi}(u^*_t)\| + \|G\|\|\nabla e_{\psi}(y^*_t)\| $.
\end{theorem}

Detailed expressions  for the constants $\kappa_1, \kappa_2,\kappa_3$ are provided in Appendix~\ref{S:completeresult}.
By comparison with the definition of $\Delta(t)$,  
Theorem~\ref{thm:finiteBasisFunctions}, $\Xi(t)$ simultaneously accounts for 
estimation error and truncation error.  We also notice that similar claims to 
those in Remark~\ref{rem:boundedness} can be derived in this case. 

\subsection{Parameter learning}
\label{sec:estimation}

We present some  methods that can be utilized in the step $\operatorname{parameter-learning}(\cdot)$ in  Framework~\ref{eq:plantModelandController}.

\emph{\textbf{Least Squares Estimator}}. Consider  the function $\phi(u)$. Suppose that at a given time $t \in \realnneg$, $K>0$ data points $(\hat{\phi}_{t_k}, u_{t_k})$ are available, with $\{t_k\}_{k = 1}^{K-1}$ the time instants where data points are received. Let $\hat{\Phi}_{t} := (\hat{\phi}_{t_1}, \ldots, \hat{\phi}_{t_K})$ be a  vector collecting the functional evaluations received up to time $t$, and $B_t \in \mathbb{R}^{K \times N}$ be a matrix with rows the regression vectors $\{b(u_{t_k})^\top, k = 1, \ldots, K\}$. Given these data points, the ordinary least squares (LS) method determines a solution to
the following optimization problems; see, e.g.,~\cite{kay1993fundamentals,hastie2009elements,Beck:2017}:
\begin{align*}
\text{ if $K\geq N$:} \quad & \hat{\alpha}_t = \arg \min_{\alpha \in \mathbb{R}^N} ~~ \norm{\hat{\Phi}_{t} - B_t \alpha}^2, &
\text{ if $K<N$:} \quad  \hat{\alpha}_t = \arg &\min_{\alpha \in \mathbb{R}^N, ~\text{s.t.:~} \hat{\Phi}_{t} = B_t\alpha} ~~ \norm{\alpha}^2 . 
\end{align*}
The above optimization problems admit a unique closed-form solution 
given by $\hat{\alpha}_t = B_t^\pinv \hat{\Phi}_{t}$, where $B_t^\pinv$ denotes the 
Moore-Penrose inverse of $B_t$. Moreover, the resulting 
approximation error admits a closed-form expression given by 
$\norm{B_t \hat{\alpha}_t -\hat{\Phi}_{t}}^2 = \norm{(I - B_t B_t^\pinv) \hat{\Phi}_{t}}^2$, which can be 
interpreted by noting that $(I - B_t B_t^\pinv)$ is the orthogonal projector
onto the null space of $B_t^\tsp$.

 A similar procedure can be utilized for the function $\psi$; in particular, we let $\hat{\Psi}_{t} := (\hat{\psi}_{t_1}, \ldots, \hat{\psi}_{t_K})$ be a  vector collecting the functional evaluations of $\psi$ received up to time $t$, and we define the matrix $D_t:= (d(y_{t_1})^\top, \dots, d(y_{t_{K}})^\top)$. Then, the LS yields the estimate $\hat{\rho}_t = D_t^\pinv \hat{\Psi}_{t}$.

\emph{\textbf{Recursive Least Squares}}. To avoid the computation of the Moore-Penrose inverse, one can utilize the recursive LS approach; we refer the reader to~\cite{ljung1999system} and~\cite{kay1993fundamentals} an overview and for the main equations of the recursive LS. See also~\cite{notarnicola2021distributed}.

\emph{\textbf{Ridge Regression}}. Given the data points $(\hat{\Psi}_{t}, B_t)$, the ridge regression involves the solution of the optimization problem $\min_{\alpha \in \mathbb{R}^N} \norm{\hat{\Phi}_{t} - B_t \alpha}^2 + \lambda_t \|\alpha\|^2$, where $\lambda_t > 0$ is a tuning parameter. While this criterion  was proposed to alleviate the singularity of $B_t^\top B_t$ when $K < N$ in~\cite{hoerl1970ridge}, the regularization $\lambda_t \|\alpha\|^2$ can be shown to impose a penalty
on the norm of $\alpha$; in fact, for a given $\lambda_t >0$, there exists $\nu_t > 0$ such that the solution of the ridge regression is equivalent to the LS with the constraint $\|\alpha\| \leq \nu_t$ as explained by, e.g.,~\cite{hastie2009elements}. The ridge regression problem admits a unique closed-form solution given by $\hat{\alpha}_t = (B_t^\top B_t + \lambda_t I)^{-1} B_t^\top \hat{\Phi}_{t}$, where $I_N$ denotes the identity matrix.  Similarly, upon receiving the $K$-th data point, the estimate of the vector $\hat{\rho}_t$ can be updated as $\hat{\rho}_t = (D_t^\top D_t + \lambda_t I)^{-1} D_t^\top \hat{\Psi}_{t}$, where $I_N$ denotes the identity matrix. To avoid the matrix inversion, a recursive strategy via the Woodbury matrix identity can be adopted.  

\emph{\textbf{Sparse Linear Regression}}. To select the basis functions that provide a parsimonious representation of 
the function, one could utilize a sparse liner regression method. This amounts to solving the problem $\min_{\alpha \in \mathbb{R}^N} \norm{\hat{\Phi}_{t} - B_t \alpha}^2 + \lambda_t \|\alpha\|_1$, where $\lambda_t > 0$ is a tuning parameter that promotes sparsity of the vector $\alpha$ as explained in~\cite{tibshirani1996regression}. The solution of the problem can be found in closed form, where the $i$-th entry of $\hat{\alpha}_t $ is given by  is $\hat{\alpha}_{i,t} = \max\{|z_{i,t}| - \lambda_t, 0 \} \textrm{sgn}(z_{i,t})$, with $z_{t} = (B_t^\top B_t)^{-1} B_t^\top \hat{\Phi}_{t}$ and where $\textrm{sgn}(\cdot)$ is the sign function.

\section{Numerical Verification}\label{S:numericalVerification}

In this section, we numerically verify the proposed Framework~\ref{eq:plantModelandController} in two cases: (i) with constant disturbance $w_t$, and (ii) with  time-varying disturbance. As an illustrative example, we utilize the LS estimator described in Section~\ref{sec:estimation}. 
We consider the cost functions $\phi(u) = \frac{1}{2} u^\top \Upsilon u +  \upsilon^\top u + r$ and $\psi(y) = \frac{1}{2} \|G u + H w_t - \xi \|^2$, 
where $ \Upsilon \in \mathbb{R}^{m \times m}$ is a symmetric positive-definite matrix, $\upsilon \in \mathbb{R}^m$ is positive element-wise, $r \in \mathbb{R}$ is positive, and $\xi \in \R^p$ is a reference output of the system. 

Since the functions are convex and quadratic, we consider the basis expansion 
in Example~\eqref{ex:parameteric} (see also~\cite{notarnicola2021distributed}); notice that to get orthonomal basis 
functions, we remove the lower-triangular part of $\Upsilon$ and we impose 
that $\Upsilon_{ij} = \Upsilon_{ji}$, for $i \neq j$. For illustration purposes, we consider a case where we estimate $\phi(u)$ (while $\psi(y)$ is known), and we consider the case $m = 4$. We generate some of the matrices of the plant and of the cost functions using normally distributed random variables; the values of the matrices are reported in Appendix~\ref{ap:simulations}. We utilize four data points $\{(u,\phi(u))\}$ as recorded data; as such, the LS is under-determined at the start of the algorithm. During the execution of the algorithm, we simulate the arrival of new data points by using a Poisson clock.  The gain of the controller satisfies the condition outlines in Theorem~\ref{thm:finiteBasisFunctions}.

Figure~2 illustrates the evolution of the error $z(t)$, defined in~\eqref{eq:trackingErrorDef}, as well as the theoretical bound provided in equation~\eqref{eq:thmResult} of Theorem~\ref{thm:finiteBasisFunctions}. Figure~2(a) illustrated the case where the disturbance is constant; new functional evaluations are received at the times marked with vertical black dotted lines. The theoretical bound depends on the estimation error, and it exhibits step changes when a new data point is received; the theoretical bound appears to be tight in the numerical simulation. 

Figure~2(a) illustrated the case where the disturbance is time varying; in this case, the bound is affected by both the estimation error and $\|\dot w_t\|$. In both cases, the numerical trajectory exhibits an exponential convergence up to an asymptotic error.   
The asymptotic error is affected by the error in the estimation of the parameters and  the variability of the disturbance.

\begin{figure}[t!]
 \centering
 \begin{minipage}[b]{0.48\linewidth}
    \centering
    \includegraphics[width=\textwidth]{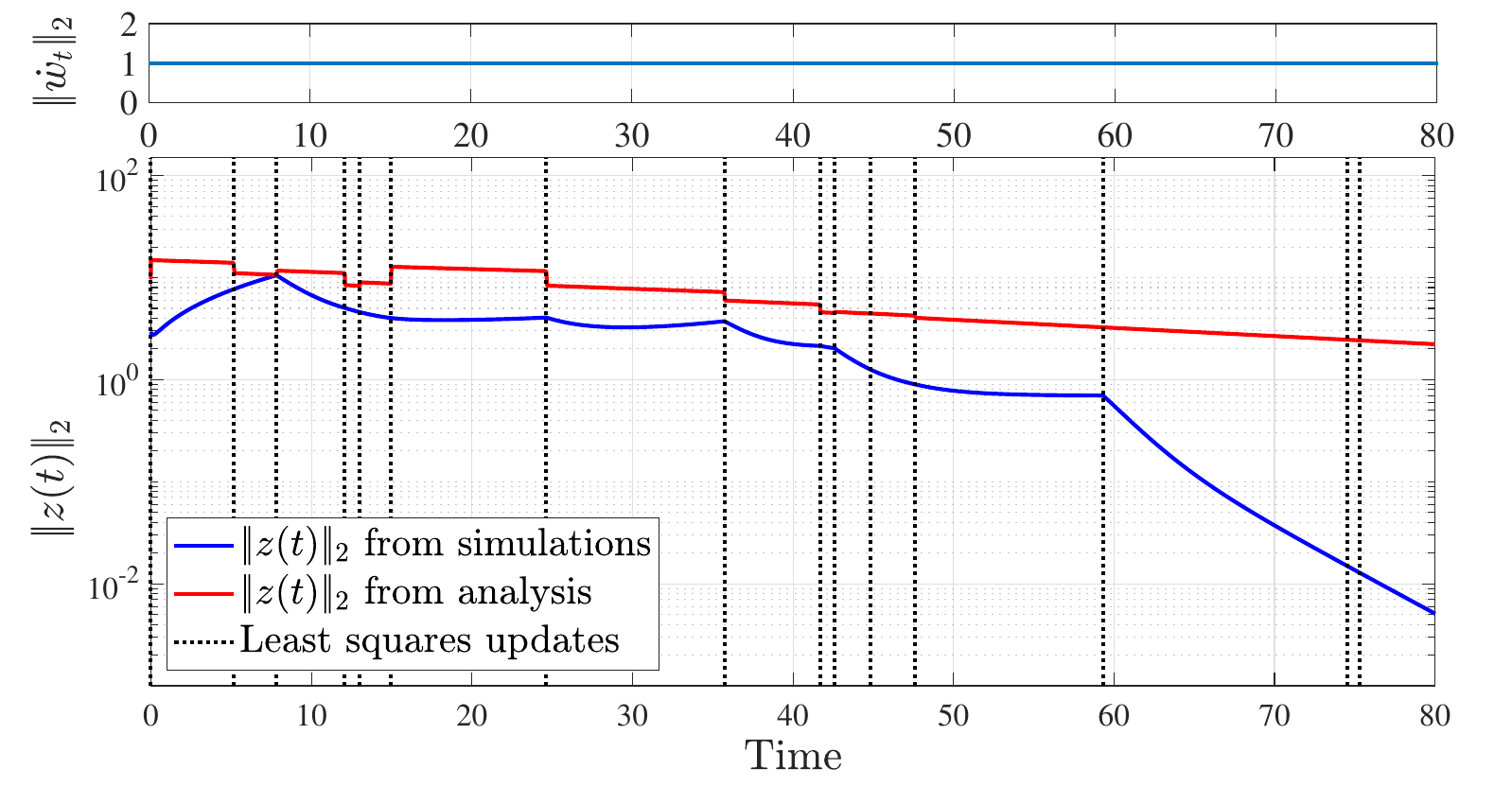}
    \small{(a)}
 \end{minipage}
 \quad
 \begin{minipage}[b]{0.48\linewidth}
    \centering
    \includegraphics[width=\textwidth]{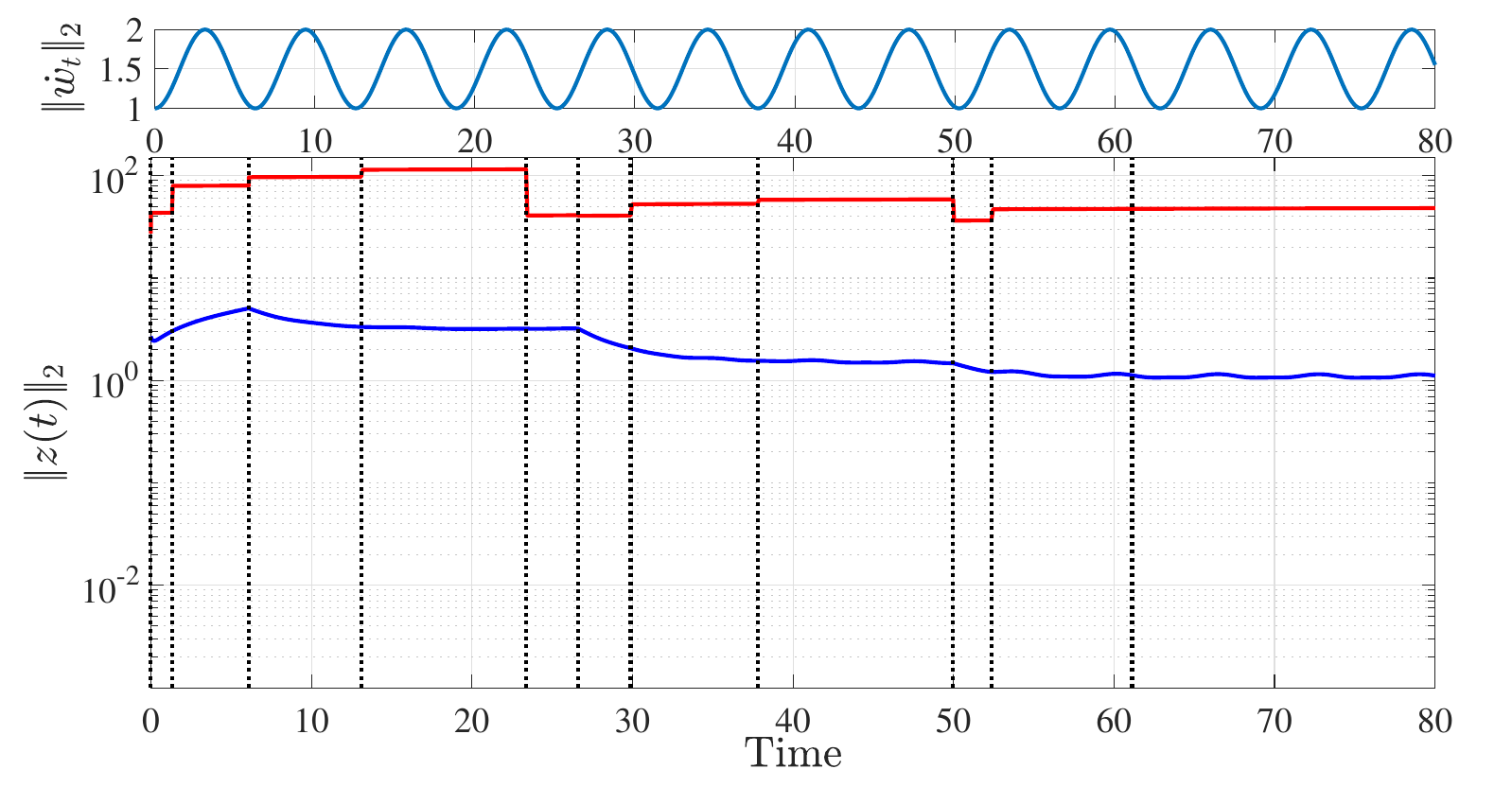} \small{(b)}
\end{minipage}
\caption{Evolution of the error $\|z(t)\|$, and theoretical bounds provided in Theorem~\ref{thm:finiteBasisFunctions}. Left: constant disturbance, $\dot w_t = 0$. Right: time-varying disturbance, $\dot w_t \neq 0$. Vertical black dotted lines represent time instants where functional evaluations are received.}
\end{figure}

\section{Conclusions}\label{S:conclusion}
We proposed a  data-enabled gradient-flow  controller to regulate an LTI dynamical system to the minimizer of an unknown functions. The controller is aided by a learning method that estimates the unknown costs from functional evaluations; to this end, appropriate  basis expansion representations (either parametric or non-parametric) are utilized. We established sufficient conditions on the estimation error and the controller gain to ensure that the error between the optimizer of the problem and the state of the closed-loop system is ultimately bounded; the error bound accounts for the functional estimation errors and the temporal variability of the unknown disturbance. Future works will look at learning methods such as concurrent learning  dynamics and neural networks.

\acks{This work was supported by the National Science Foundation (NSF) through the Awards CMMI 2044946 and CAREER 1941896. }

\bibliography{alias,bibliography}

\begin{thebibliography}{53}
\providecommand{\natexlab}[1]{#1}
\providecommand{\url}[1]{\texttt{#1}}
\expandafter\ifx\csname urlstyle\endcsname\relax
  \providecommand{\doi}[1]{doi: #1}\else
  \providecommand{\doi}{doi: \begingroup \urlstyle{rm}\Url}\fi

\bibitem[Angeli et~al.(2003)Angeli, Sontag, and Wang]{angeli2003input}
D.~Angeli, E.~D. Sontag, and Y.~Wang.
\newblock Input-to-state stability with respect to inputs and their
  derivatives.
\newblock \emph{International Journal of Robust and Nonlinear Control:
  IFAC-Affiliated Journal}, 13\penalty0 (11):\penalty0 1035--1056, 2003.

\bibitem[Ariyur and Krsti\'c(2003)]{ariyur2003real}
K.~B. Ariyur and M.~Krsti\'c.
\newblock \emph{Real-time optimization by extremum-seeking control}.
\newblock John Wiley \& Sons, 2003.

\bibitem[Bazerque and Giannakis(2013)]{bazerque2013nonparametric}
J.~Bazerque and G.~B. Giannakis.
\newblock Nonparametric basis pursuit via sparse kernel-based learning: A
  unifying view with advances in blind methods.
\newblock \emph{IEEE Signal Processing Magazine}, 30\penalty0 (4):\penalty0
  112--125, 2013.

\bibitem[Beck(2014)]{Beck:2014}
A.~Beck.
\newblock \emph{Introduction to Nonlinear Optimization: Theory, Algorithms, and
  Applications with MATLAB}.
\newblock MOS-SIAM Series on Optimization, Technion-Israel Institute of
  Technology, Kfar Saba, Israel, 2014.

\bibitem[Beck(2017)]{Beck:2017}
A.~Beck.
\newblock \emph{First Order Methods in Optimization}.
\newblock MOS-SIAM Series on Optimization, Tel-Aviv University, Tel-Aviv,
  Israel, 2017.

\bibitem[Bianchin et~al.(2020)Bianchin, Poveda, and
  Dall'Anese]{bianchin2020online}
G.~Bianchin, J.~I. Poveda, and E.~Dall'Anese.
\newblock Online optimization of switched {LTI} systems using continuous-time
  and hybrid accelerated gradient flows.
\newblock \emph{arXiv preprint arXiv:2008.03903}, 2020.

\bibitem[Bianchin et~al.(2021{\natexlab{a}})Bianchin, Cort\'es, Poveda, and
  Dall'Anese]{bianchin2021time}
G.~Bianchin, J.~Cort\'es, J.~I. Poveda, and E.~Dall'Anese.
\newblock Time-varying optimization of {LTI} systems via projected primal-dual
  gradient flows.
\newblock \emph{IEEE Trans. on Control of Networked Systems},
  2021{\natexlab{a}}.

\bibitem[Bianchin et~al.(2021{\natexlab{b}})Bianchin, Dall'Anese, Poveda,
  Jacobson, Carlton, and Buchwald]{bianchin2021planning}
G.~Bianchin, E.~Dall'Anese, J.~I. Poveda, D.~Jacobson, E.~J. Carlton, and A.~G.
  Buchwald.
\newblock Planning a return to normal after the {COVID}-19 pandemic:
  Identifying safe contact levels via online optimization.
\newblock \emph{arXiv preprint arXiv:2109.06025}, 2021{\natexlab{b}}.

\bibitem[Bourgin et~al.(2019)Bourgin, Peterson, Reichman, Russell, and
  Griffiths]{bourgin2019cognitive}
D.~D. Bourgin, J.~C. Peterson, D.~Reichman, S.~J. Russell, and T.~L. Griffiths.
\newblock Cognitive model priors for predicting human decisions.
\newblock In \emph{International Conference on Machine Learning}, pages
  5133--5141, 2019.

\bibitem[Brunner et~al.(2012)Brunner, D{\"u}rr, and
  Ebenbauer]{brunner2012feedback}
F.~D. Brunner, H.-B. D{\"u}rr, and C.~Ebenbauer.
\newblock Feedback design for multi-agent systems: A saddle point approach.
\newblock In \emph{{IEEE} Conf.\ on Decision and Control}, pages 3783--3789,
  2012.

\bibitem[Chowdhary and Johnson(2010)]{chowdhary2010concurrent}
G.~Chowdhary and E.~Johnson.
\newblock Concurrent learning for convergence in adaptive control without
  persistency of excitation.
\newblock In \emph{{IEEE} Conf.\ on Decision and Control}, pages 3674--3679.
  IEEE, 2010.

\bibitem[Chowdhary et~al.(2013)Chowdhary, Yucelen, M{\"u}hlegg, and
  Johnson]{chowdhary2013concurrent}
G.~Chowdhary, T.~Yucelen, M.~M{\"u}hlegg, and E.~N. Johnson.
\newblock Concurrent learning adaptive control of linear systems with
  exponentially convergent bounds.
\newblock \emph{International Journal of Adaptive Control and Signal
  Processing}, 27\penalty0 (4):\penalty0 280--301, 2013.

\bibitem[Colombino et~al.(2020)Colombino, Dall’Anese, and
  Bernstein]{MC-ED-AB:20}
M.~Colombino, E.~Dall’Anese, and A.~Bernstein.
\newblock Online optimization as a feedback controller: Stability and tracking.
\newblock \emph{IEEE Trans. on Control of Network Systems}, 7\penalty0
  (1):\penalty0 422--432, 2020.

\bibitem[Cothren et~al.(2021)Cothren, Bianchin, and
  Dall'Anese]{cothren2021Data}
L.~Cothren, G.~Bianchin, and E~Dall'Anese.
\newblock Data-enabled gradient flow as feedback controller: Regulation of
  linear dynamical systems to minimizers of unknown functions (extended
  version).
\newblock \emph{arXiv preprint}, 2021.
\newblock https://arxiv.org/abs/2112.01652.

\bibitem[Davison(1976)]{ED:76}
E.~Davison.
\newblock The robust control of a servomechanism problem for linear
  time-invariant multivariable systems.
\newblock \emph{IEEE Trans. on Automatic Control}, 21\penalty0 (1):\penalty0
  25--34, 1976.

\bibitem[Devraj and Meyn(2017)]{devraj2017zap}
A.~M Devraj and S.~P Meyn.
\newblock Zap {Q}-learning.
\newblock In \emph{Proceedings of the 31st International Conference on Neural
  Information Processing Systems}, pages 2232--2241, 2017.

\bibitem[Fabiani et~al.(2021)Fabiani, Simonetto, and
  Goulart]{fabiani2021learning}
F.~Fabiani, A.~Simonetto, and P.~J. Goulart.
\newblock Learning equilibria with personalized incentives in a class of
  nonmonotone games.
\newblock \emph{arXiv preprint arXiv:2111.03854}, 2021.

\bibitem[G{\"u}ndling et~al.(2020)G{\"u}ndling, Hopp, and
  Weihe]{gundling2020efficient}
F.~G{\"u}ndling, F.~Hopp, and K.~Weihe.
\newblock Efficient monitoring of public transport journeys.
\newblock \emph{Public Transport}, 12\penalty0 (3):\penalty0 631--645, 2020.

\bibitem[Hastie et~al.(2009)Hastie, Tibshirani, and
  Friedman]{hastie2009elements}
T.~Hastie, R.~Tibshirani, and J.~Friedman.
\newblock \emph{The elements of statistical learning}.
\newblock Springer, 2009.

\bibitem[Hauswirth et~al.(2021{\natexlab{a}})Hauswirth, Bolognani, Hug, and
  D{\"o}rfler]{hauswirth2020timescale}
A.~Hauswirth, S.~Bolognani, G.~Hug, and F.~D{\"o}rfler.
\newblock Timescale separation in autonomous optimization.
\newblock \emph{IEEE Trans. on Automatic Control}, 66\penalty0 (2):\penalty0
  611--624, 2021{\natexlab{a}}.

\bibitem[Hauswirth et~al.(2021{\natexlab{b}})Hauswirth, Bolognani, Hug, and
  D{\"o}rfler]{hauswirth2021optimization}
A.~Hauswirth, S.~Bolognani, G.~Hug, and F.~D{\"o}rfler.
\newblock Optimization algorithms as robust feedback controllers.
\newblock \emph{arXiv preprint arXiv:2103.11329}, 2021{\natexlab{b}}.

\bibitem[Hirata et~al.(2014)Hirata, Hespanha, and Uchida]{Hirata}
K.~Hirata, J.~P. Hespanha, and K.~Uchida.
\newblock Real-time pricing leading to optimal operation under distributed
  decision makings.
\newblock In \emph{Proc. of American Control Conf.}, Portland, OR, June 2014.

\bibitem[Hoerl and Kennard(1970)]{hoerl1970ridge}
A.~E. Hoerl and R.~W. Kennard.
\newblock Ridge regression: Biased estimation for nonorthogonal problems.
\newblock \emph{Technometrics}, 12\penalty0 (1):\penalty0 55--67, 1970.

\bibitem[Hornik et~al.(1989)Hornik, Stinchcombe, and
  White]{hornik1989multilayer}
K.~Hornik, M.~Stinchcombe, and H.~White.
\newblock Multilayer feedforward networks are universal approximators.
\newblock \emph{Neural networks}, 2\penalty0 (5):\penalty0 359--366, 1989.

\bibitem[Jin and Lavaei(2020)]{MJ-JL:20}
M.~Jin and J.~Lavaei.
\newblock Stability-certified reinforcement learning: A control-theoretic
  perspective.
\newblock \emph{IEEE Access}, 8:\penalty0 229086--229100, 2020.

\bibitem[Jokic et~al.(2009)Jokic, Lazar, and Van
  Den~Bosch]{Jokic2009controller}
A.~Jokic, M.~Lazar, and P.~P.-J. Van Den~Bosch.
\newblock On constrained steady-state regulation: Dynamic {KKT} controllers.
\newblock \emph{IEEE Trans. on Automatic Control}, 54\penalty0 (9):\penalty0
  2250--2254, 2009.

\bibitem[Karg and Lucia(2020)]{karg2020stability}
B.~Karg and S.~Lucia.
\newblock Stability and feasibility of neural network-based controllers via
  output range analysis.
\newblock In \emph{IEEE Conference on Decision and Control}, pages 4947--4954,
  2020.

\bibitem[Karimi et~al.(2016)Karimi, Nutini, and Schmidt]{karimi2016linear}
H.~Karimi, J.~Nutini, and M.~Schmidt.
\newblock Linear convergence of gradient and proximal-gradient methods under
  the {Polyak-{\L}ojasiewicz} condition.
\newblock In \emph{Joint European Conference on Machine Learning and Knowledge
  Discovery in Databases}, pages 795--811. Springer, 2016.

\bibitem[Kay(1993)]{kay1993fundamentals}
S.~M. Kay.
\newblock \emph{Fundamentals of statistical signal processing: estimation
  theory}.
\newblock Prentice-Hall, Inc., 1993.

\bibitem[Khalil(2002)]{Khalil:1173048}
H.~K. Khalil.
\newblock \emph{Nonlinear Systems; 3rd ed.}
\newblock Prentice-Hall, Upper Saddle River, NJ, 2002.

\bibitem[Krstic and Wang(2000)]{krstic2000stability}
M.~Krstic and H.-H. Wang.
\newblock Stability of extremum seeking feedback for general nonlinear dynamic
  systems.
\newblock \emph{Automatica}, 36\penalty0 (4):\penalty0 595--602, 2000.

\bibitem[Lawrence et~al.(2018)Lawrence, Nelson, Mallada, and
  Simpson-Porco]{lawrence2018optimal}
L.~S.~P. Lawrence, Z.~E. Nelson, E.~Mallada, and J.~W. Simpson-Porco.
\newblock Optimal steady-state control for linear time-invariant systems.
\newblock In \emph{{IEEE} Conf.\ on Decision and Control}, pages 3251--3257,
  December 2018.

\bibitem[Lawrence et~al.(2021)Lawrence, Simpson-Porco, and
  Mallada]{lawrence2018linear}
L.~S.~P. Lawrence, J.~W. Simpson-Porco, and E.~Mallada.
\newblock Linear-convex optimal steady-state control.
\newblock \emph{IEEE Trans. on Automatic Control}, 2021.
\newblock ({To} appear).

\bibitem[Ljung(1999)]{ljung1999system}
L.~Ljung.
\newblock \emph{System identification: theory for the user}.
\newblock 2nd edition Prentice-Hall, Upper Saddle River, NJ, 1999.

\bibitem[Luan et~al.(2013)Luan, Sharp, and LaRoy]{luan2013data}
W.~Luan, D.~Sharp, and S.~LaRoy.
\newblock Data traffic analysis of utility smart metering network.
\newblock In \emph{IEEE Power \& Energy Society General Meeting}, pages 1--4,
  2013.

\bibitem[Luo et~al.(2020)Luo, Zhang, and Zavlanos]{luo2020socially}
X.~Luo, Y.~Zhang, and M.~M. Zavlanos.
\newblock Socially-aware robot planning via bandit human feedback.
\newblock In \emph{International Conf. on Cyber-Physical Systems}, pages
  216--225. IEEE, 2020.

\bibitem[Marchi et~al.(2022)Marchi, Bunton, Gharesifard, and
  Tabuada]{Marchi2022}
M.~Marchi, J.~Bunton, B.~Gharesifard, and P.~Tabuada.
\newblock Safety and stability guarantees for control loops with deep learning
  perception.
\newblock \emph{IEEE Control Systems Letters}, 6:\penalty0 1286--1291, 2022.

\bibitem[Menta et~al.(2018)Menta, Hauswirth, Bolognani, Hug, and
  D{\"o}rfler]{menta2018stability}
S.~Menta, A.~Hauswirth, S.~Bolognani, G.~Hug, and F.~D{\"o}rfler.
\newblock Stability of dynamic feedback optimization with applications to power
  systems.
\newblock In \emph{Annual Conf. on Communication, Control, and Computing},
  pages 136--143, 2018.

\bibitem[Munir et~al.(2013)Munir, Stankovic, Liang, and Lin]{munir2013cyber}
S.~Munir, J.~A. Stankovic, C.-J.~M. Liang, and S.~Lin.
\newblock Cyber physical system challenges for human-in-the-loop control.
\newblock In \emph{8th International Workshop on Feedback Computing}, 2013.

\bibitem[Notarnicola et~al.(2021)Notarnicola, Simonetto, Farina, and
  Notarstefano]{notarnicola2021distributed}
I.~Notarnicola, A.~Simonetto, F.~Farina, and G.~Notarstefano.
\newblock Distributed personalized gradient tracking with convex parametric
  models, 2021.

\bibitem[Ospina et~al.(2020)Ospina, Simonetto, and
  Dall’Anese]{ospina2020personalized}
A.~M. Ospina, A.~Simonetto, and E.~Dall’Anese.
\newblock Personalized demand response via shape-constrained online learning.
\newblock In \emph{IEEE International Conference on Communications, Control,
  and Computing Technologies for Smart Grids}, pages 1--6, 2020.

\bibitem[Poveda et~al.(2021)Poveda, Benosman, and Vamvoudakis]{poveda2021data}
J.~I. Poveda, M.~Benosman, and K.~G. Vamvoudakis.
\newblock Data-enabled extremum seeking: a cooperative concurrent
  learning-based approach.
\newblock \emph{International Journal of Adaptive Control and Signal
  Processing}, 35\penalty0 (7):\penalty0 1256--1284, 2021.

\bibitem[Robey et~al.()Robey, H., Lindemann, Zhang, Dimarogonas, Tu, and
  Matni]{robey2020learning}
A.~Robey, Hu~H., L.~Lindemann, H.~Zhang, D.~Dimarogonas, S.~Tu, and N.~Matni.
\newblock Learning control barrier functions from expert demonstrations.
\newblock In \emph{{IEEE} Conf.\ on Decision and Control}.

\bibitem[Simonetto et~al.(2021)Simonetto, Dall’Anese, Monteil, and
  Bernstein]{simonetto2021personalized}
A.~Simonetto, E.~Dall’Anese, J.~Monteil, and A.~Bernstein.
\newblock Personalized optimization with user’s feedback.
\newblock \emph{Automatica}, 131:\penalty0 109767, 2021.

\bibitem[Sontag(2022)]{sontag2021remarks}
E.~D. Sontag.
\newblock Remarks on input to state stability on open sets, motivated by
  perturbed gradient flows in model-free learning.
\newblock \emph{Systems \& Control Letters}, 161:\penalty0 105138, 2022.

\bibitem[Sontag and Wang(1997)]{sontag1997output}
E.~D. Sontag and Y.~Wang.
\newblock Output-to-state stability and detectability of nonlinear systems.
\newblock \emph{Systems \& Control Letters}, 29\penalty0 (5):\penalty0
  279--290, 1997.

\bibitem[Tabuada and Gharesifard(2020)]{tabuada2020universal}
P.~Tabuada and B.~Gharesifard.
\newblock Universal approximation power of deep residual neural networks via
  nonlinear control theory.
\newblock \emph{arXiv preprint arXiv:2007.06007}, 2020.

\bibitem[Taylor et~al.(2020)Taylor, Singletary, Yue, and
  Ames]{taylor2020learning}
A.~Taylor, A.~Singletary, Y.~Yue, and A.~Ames.
\newblock Learning for safety-critical control with control barrier functions.
\newblock In \emph{Learning for Dynamics and Control}, pages 708--717. PMLR,
  2020.

\bibitem[Teel and Popovic(2001)]{teel2001solving}
A.~R. Teel and D.~Popovic.
\newblock Solving smooth and nonsmooth multivariable extremum seeking problems
  by the methods of nonlinear programming.
\newblock In \emph{American Control Conference}, volume~3, pages 2394--2399,
  2001.

\bibitem[Tibshirani(1996)]{tibshirani1996regression}
R.~Tibshirani.
\newblock Regression shrinkage and selection via the lasso.
\newblock \emph{Journal of the Royal Statistical Society}, 58\penalty0
  (1):\penalty0 267--288, 1996.

\bibitem[van Nunen et~al.(2017)van Nunen, Esposto, Saberi, and
  Paardekooper]{Nunen}
E.~van Nunen, F.~Esposto, A.~K. Saberi, and J.-P. Paardekooper.
\newblock Evaluation of safety indicators for truck platooning.
\newblock In \emph{IEEE Intelligent Vehicles Symposium}, pages 1013--1018,
  2017.

\bibitem[Yin et~al.(2021)Yin, Seiler, and Arcak]{HY-PS-MA:21}
H.~Yin, P.~Seiler, and M.~Arcak.
\newblock Stability analysis using quadratic constraints for systems with
  neural network controllers.
\newblock \emph{IEEE Trans. on Automatic Control}, 2021.
\newblock (Early Access).

\bibitem[Zheng et~al.(2020)Zheng, Simpson-Porco, and
  Mallada]{zheng2019implicit}
T.~Zheng, J.~Simpson-Porco, and E.~Mallada.
\newblock Implicit trajectory planning for feedback linearizable systems: A
  time-varying optimization approach.
\newblock In \emph{American control Conference}, pages 4677--4682, 2020.

\end{thebibliography}

\section{Complete Results for Theorems~\ref{thm:finiteBasisFunctions} and~\ref{thm:truncation}}
\label{S:completeresult}

Here we provide more details on the Theorems~\ref{thm:finiteBasisFunctions} and~\ref{thm:truncation}. In the main results, we define the constants $\kappa_1, \kappa_2$, $\kappa_3 > 0$ as: 
\begin{align}
    \kappa_1 &:= \sqrt{\frac{c_1}{c_2}} \\ \kappa_2 & := \frac{c_4}{2 \sqrt{c_1}} \\
    \kappa_3 & := \frac{c_5}{2\sqrt{c_1}}
 \end{align}   
where: 
\begin{align}
        c_1 &:= \min \left\{\frac{(1 - \theta)}{\eta}\frac{\mu_{u}}{2},\frac{\theta}{\eta}\underline{\lambda}\left(P\right)\right\},  \\
        c_2 & := \max \left\{ \frac{(1 - \theta)}{\eta}\frac{\ell}{2}, \frac{\theta}{\eta} \bar{\lambda} (P) \right\}, \\
        c_4 &:= \sqrt{\eta} \max \left\{\ell \sqrt{\frac{2}{\mu_u}},  \frac{2\|PA^{-1}B\|}{\sqrt{\underline{\lambda}(P)}}\right\}, \\ 
        c_5 & := \frac{2 \|P^\top A^{-1} E\|}{\sqrt{\eta} \sqrt{\underline{\lambda}(P)}}.
\end{align}

\section{Proofs of the Main Results}
\label{S:proofs}

To prove our main results, we combine arguments from perturbation theory and singular perturbation theory (respectively, \cite[Ch. 9, Sec. 9.3]{Khalil:1173048} and \cite[Ch. 11, S. 11.5]{Khalil:1173048}). First, we 
derive a Lyapunov function inspired from \cite[Thm. 11.3, Lem. 9.3]{Khalil:1173048} and we characterize the choices of the controller gain $\eta$ that guarantee stability of the controlled system in the absence of 
disturbances.
Second, we use 
Lemma 9.5 \cite[Sec. 9.3]{Khalil:1173048} to characterize the transient behavior
behavior of the controller error in the presence of disturbances.

\subsection{Notation for Cost Functions}

Given \( \phi(u) = \sum_{i=1}^{N} \alpha_{i} b_i(u) + \sum_{i=N+1}^{N_b} b_i(u) \alpha_i\), $N < \infty$, with basis functions $b_i: \mathbb{R}^m \rightarrow \mathbb{R}$, denote the column vectors \( \alpha := [\alpha_{1}, \alpha_{2}, \dots, \alpha_{N}]^\top \), \( b(u) := [b_1(u), b_2(u), \dots, b_{N}(u)]^\top \), and $e_{\phi}(u) := \sum_{i=N+1}^{N_b} b_i(u) \alpha_i$. Then, under Assumption \ref{as:Phi} the function 
$\phi(\cdot)$ can be rewritten as \[\phi(u) = b(u)^\top \alpha + e_\phi(u). \]

The gradient of $\phi(u)$ is given by $\nabla \phi(u) = \sum_{i = 1}^{N} \alpha_i \nabla b_i (u) + \nabla e_\phi(u)$, and the estimated gradient based on an estimate $\hat \alpha_t$ of $\alpha$ available at time $t$ is $\nabla \hat{\phi}(u) = \sum_{i = 1}^{N} \hat{\alpha}_{i,t} \nabla b_i (u)$. Let $\nabla b(u)$ denote the Jacobian of $b(u)$, or:
\begin{equation*}
\nabla b := \left[
    \begin{array}{cccc}
    \frac{\partial b_1}{\partial u_1} & \ldots & \ldots &  \frac{\partial b_1}{\partial u_m} \\
    \vdots & & & \vdots \\
    \frac{\partial b_{N}}{\partial u_1} & \ldots & \ldots &  \frac{\partial b_{N}}{\partial u_m} 
    \end{array}
    \right]  .
\end{equation*}
Then, we write the gradient of the true function $\phi$ and the estimated function $\hat{\phi}$, respectively, as, \[\nabla \phi(u) = \nabla b(u)^\top \alpha + \nabla e_\phi(u) \text{ and } \nabla \hat{\phi}(u) = \nabla b(u)^\top \hat{\alpha}_t.\]

Similarly, given \( \psi(y) = \sum_{j=1}^{M} \rho_{j} d_j(y) + \sum_{j=M+1}^{M_b}  d_j(y) \rho_j\), $M < \infty$, with basis functions $d_j: \mathbb{R}^p \rightarrow \mathbb{R}$, 
denote the column vectors \( \rho := [\rho_{1}, \rho_{2}, \dots, \rho_{M}]^\top \), \( d(y) := [d_1(y), d_2(y), \dots, d_{M}(y)]^\top \), and truncation error $e_{\psi}(y) :=  \sum_{j=M+1}^{M_b}  d_j(y) \rho_j$ to rewrite the cost function as
\[\psi(y) = d(y)^\top \rho + e_\psi(y). \]

The gradient of $\psi(y)$ is given by $\nabla \psi(y) = \sum_{i = 1}^{M} \rho_i \nabla d_i (y) + \nabla e_\psi(y)$, and the estimated gradient based on an estimate $\hat \rho_t$ of $\rho$ available at time $t$ is $\nabla \hat{\psi}(y) = \sum_{i = 1}^{M} \hat{\rho}_{i,t} \nabla d_i (y)$. Let $\nabla d(u)$ denote the Jacobian of $d(y)$, namely,
\begin{equation*}
\nabla d := \left[
    \begin{array}{cccc}
    \frac{\partial d_1}{\partial y_1} & \ldots & \ldots &  \frac{\partial d_1}{\partial y_p} \\
    \vdots & & & \vdots \\
    \frac{\partial d_{M}}{\partial y_1} & \ldots & \ldots &  \frac{\partial d_{M}}{\partial y_p} 
    \end{array}
    \right]  .
\end{equation*}
Then, we write the gradient of the true function $\psi$ and the estimated function $\hat{\psi}$, respectively, as, \[\nabla \psi(y) = \nabla d(y)^\top \rho +\nabla e_\psi(y) \text{ and } \nabla \hat{\psi}(y) = \nabla d(y)^\top \hat{\rho}_t.\]

We recall that we consider two cases. In \emph{\textbf{(Case~1)}} the functions $\phi$ and $\psi$ are both represented by a finite number of basis functions, and we set $N = N_b$ and $M = M_b$ basis functions;  in this case, $e_\phi(u) = 0$ and $e_\psi(y) = 0$. In  \emph{\textbf{(Case~2)}}, we have $N \ll N_b$ and $M \ll M_b$ (where $N_b$ and $M_b$ may be large or even $\infty$). In this section, we outline the main proof for case \emph{\textbf{(Case~2)}}; the proof for \emph{\textbf{(Case~1)}} follows directly, and we will specify relevant modifications whenever needed.   

\subsection{Perturbed Gradient Flow and Singular-Perturbation Model}

Consider the controller, which is based on a perturbed gradient flow: 
\begin{equation}
\label{eq:controller_ap}
\dot u = -\eta \left(\nabla \hat \phi(u) + \nabla \hat \psi(y)\right) = -\eta \left( \nabla b(u)^\top \hat \alpha_t + \nabla d(y)^\top \hat \rho_t \right).
\end{equation}
Rewrite~\eqref{eq:controller_ap} in terms of a nominal term (i.e., the true gradient) and an error term (i.e., the perturbation) by first adding and subtracting the true gradients, \( \nabla \phi(u) =  \nabla b(u)^\top \alpha + \nabla e_\phi(u) \) and \( \nabla \psi(y) = \nabla d(y)^\top \rho +\nabla e_\psi(y)\): 
\begin{equation*}
\dot u = -\eta \left(\nabla \phi(u) + \nabla \psi(y)\right) + \eta \left(\nabla b(u)^\top(\alpha - \hat \alpha_t) + \nabla d(y)^\top(\rho - \hat \rho_t) + \nabla e_\phi(u) + \nabla_\psi(y)\right).
\end{equation*}
Now, add and subtract \( \eta \nabla b(u^*)^\top(\alpha - \hat \alpha_t) + \eta \nabla e_\phi(u^*) \) and \( \eta \nabla d(y^*)^\top(\rho - \hat \rho_t) + \eta \nabla e_\psi(y^*)\) to obtain:
\begin{align}
\dot u = -\eta & \left\{\nabla \phi(u) + \nabla \psi(y)\right\}  + \eta \{\left(\nabla b(u) - \nabla b(u^*)\right)^\top(\alpha - \hat \alpha_t) + (\nabla d(y) - \nabla d(y^*))^\top (\rho - \hat \rho_t) \nonumber \\
& + (\nabla e_\phi(u) - \nabla e_\phi(u^*)) + (\nabla e_\psi(y) - \nabla e_\psi(y^*))  \nonumber \\
& + \nabla b(u^*)^\top(\alpha - \hat \alpha_t) + \nabla d(y^*)^\top(\rho - \hat \rho_t)  + \nabla e_\phi(u^*) + \nabla e_\psi(y^*) \}. \label{eq:nominalperturbed}
\end{align}
This representation is in line with the model of~\cite[Ch. 9, Ch. 11]{Khalil:1173048}, in which the use of a nominal and error term is inspired by Ch. 9, and the use of controller gain, $\eta$, follows singular perturbation theory. With this representation, the interconnection between plant and controller can be rewritten as: 
\begin{equation}
\begin{split}
\label{eq:plantAndController}
\dot x &=  A x + B u + E w_t, \\
y &= C x + D w_t,\\
\dot u &= -\eta \left\{\nabla \phi(u) + \nabla \psi(y)\right\} + \eta \{\left(\nabla b(u) - \nabla b(u^*)\right)^\top(\alpha - \hat \alpha_t) + (\nabla d(y) - \nabla d(y^*))^\top (\rho - \hat \rho_t) \\
&~~~~~~~~~~~~~~~~~ + (\nabla e_\phi(u) - \nabla e_\phi(u^*)) + (\nabla e_\psi(y) - \nabla e_\psi(y^*))  \\
&~~~~~~~~~~~~~~~~~ + \nabla b(u^*)^\top(\alpha - \hat \alpha_t) + \nabla d(y^*)^\top(\rho - \hat \rho_t)  + \nabla e_\phi(u^*) + \nabla e_\psi(y^*) \}.
\end{split}
\end{equation}

\subsection{Change of Variables}

We begin by performing a change of variables to shift the equilibrium point of
the plant  to the origin. For any $w_t$, let $x_{eq}$ be the (unique) 
equilibrium of \eqref{eq:plantModel}, and define 
$x' := x - x_{eq} = x - (-A^{-1}B u - A^{-1}E w_t)$ to obtain 
\begin{align*}
    \dot x = \dot{x}' + \dot{x}_{eq} = A(x'+ x_{eq}) + B u + E w_t \implies \dot x' + \dot x_{eq} = Ax' + Ax_{eq} + B u + E w_t,
\end{align*}    
and    
\begin{align*}
    \dot x' + \frac{\partial x_{eq}}{\partial u}\dot u + \frac{\partial x_{eq}}{\partial w_t}\dot w_t = Ax' + (- Bu - E w_t) + Bu + Ew_t
    \implies \dot x' = Ax + A^{-1}B \dot u + A^{-1}E \dot w_t.
\end{align*}
The system $\dot x' = Ax + A^{-1}B \dot u + A^{-1}E \dot w_t$ is the so-called ``boundary-layer system;'' see~\cite[Chapter 11]{Khalil:1173048}. 

Next, we rewrite the controller dynamics in the new variables.
Recall that at equilibrium, $y = Gu + H w_t $, with $G = -CA^{-1}B$ and $H = D - CA^{-1}E$.
Also, we know that $y = Cx + Dw_t.$ Use these facts for the following change of variables: 
\begin{align*}
    \psi(y) &= \psi(Cx + D w_t) = \psi(C(x' + x_{eq}) + D w_t)\\
    &= \psi(C x' -CA^{-1}Bu - CA^{-1}E w_t + D w_t)\\
    &= \psi( C x' + G u + H w_t).
\end{align*}
Then, the gradient can be written as $G^\top \nabla \psi(C x' + G u + H w_t)$. 
For brevity, hereafter we write $y$ to mean $y = C x' + G u + H w_t$ and $y^* = C x' + G u^* + H w_t$ for the error terms in the controller 
With these changes of variables, the controller can be written as:
\begin{align*}
    \dot u &= - \eta (\nabla \phi(u) + G^\top \nabla \psi (Cx' + G u - H w_t)) \\
    &~~~~~~~~~~~~~~~~~ + \eta \{\left(\nabla b(u) - \nabla b(u^*)\right)^\top(\alpha - \hat \alpha_t) + (G^\top \nabla d(y) - G^\top \nabla d(y^*))^\top (\rho - \hat \rho_t) \\
&~~~~~~~~~~~~~~~~~ + (\nabla e_\phi(u) - \nabla e_\phi(u^*)) + (G^\top \nabla e_\psi(y) - G^\top \nabla e_\psi(y^*))  \\
&~~~~~~~~~~~~~~~~~ + \nabla b(u^*)^\top(\alpha - \hat \alpha_t) + \nabla d(y^*)^\top(\rho - \hat \rho_t)  + \nabla e_\phi(u^*) + G^\top \nabla e_\psi(y^*) \}
\end{align*}
where we have again a nominal gradient flow and the perturbation. This system is the ``reduced system'' in our singular-perturbation setup; see~ \cite[Chapter 11]{Khalil:1173048}. 

In what follows, we denote in compact form
\begin{align*}
m(u, x) &:= -\nabla \phi(u) - G^\top \nabla \psi(Cx +Dw_t), \\
n(u, x) &:= A x + B u + E w_t, \\
e(u,y) &:= (\nabla b(u) - \nabla b(u^*))^\top (\alpha - \hat \alpha_t) + (G^\top \nabla d(y) -  G^\top\nabla d(y^*))^\top(\rho - \hat \rho_t) \\
        & ~~ + (\nabla e_\phi(u) + \nabla e_\phi(u^*)) + (G^\top \nabla e_\psi(y) - G^\top \nabla e_\psi(y^*))\\
        & ~~ + \nabla b(u^*)^\top (\alpha - \hat \alpha_t) + G^\top \nabla d(y^*)^\top (\rho + \hat \rho_t) + \nabla e_\phi(u^*) + G^\top \nabla e_\psi(y^*).
\end{align*}
Accordingly, the dynamics in the new variables read as:
\begin{align*}
\dot x^\prime &= n(u, x' + x_{eq}) - \eta \frac{\partial x_{eq}}{\partial u} \dot u - \frac{\partial x_{eq}}{\partial w_t} \dot w_t, &
\dot u &= \eta \, m(u, x' + x_{eq}) + \eta \, e(u,y).
\end{align*}

\subsection{Lyapunov Functions and Bounds}
For the boundary-layer system, propose the following Lyapunov function:  
\begin{equation}
\label{eq:BoundaryLayerLyapunov}
W(x') = x'^\top P x' 
\end{equation}
where, for any positive definite matrix $Q \in R^{n \times n}, $ $P \succ 0$ is the solution to the Lyapunov equation $A^\top P + P A = -Q$ as in Assumption \ref{ass:stabilityPlant}. Notice that \eqref{eq:BoundaryLayerLyapunov} satisfies the following quadratic bounds: 
\begin{align}
\label{eq:quadraticBoundW}
    \underline{\lambda}(P) \|x'\|^2 \leq W(x') \leq \bar \lambda(P) \|x'\|^2.
\end{align}
For the reduced system, we utilize the Lyapunov function:  
\begin{equation}
\label{eq:ReducedLyapunov}
V(t, u; w) = \phi(u) - \phi(u^*) + \psi(G u + H w_t) - \psi(G u^* + H w_t).
\end{equation}
Note that \eqref{eq:ReducedLyapunov} satisfies the bounds:  
\begin{align}
\label{eq:quadraticBoundV}
    \frac{\mu_u}{2}\|u - u^*\|^2 \leq V(u) \leq \frac{\ell_u + \ell_y\|G\|^2}{2}\|u - u^*\|^2.
\end{align}
then, combining \eqref{eq:BoundaryLayerLyapunov} and \eqref{eq:ReducedLyapunov}, we obtain the composite Lyapunov function
\begin{equation}
\label{eq:compositeLyapunov}
\nu (t, u, x'; w) = (1 - \theta) \frac{1}{\eta}V(t, u; w) + \theta \frac{1}{\eta} W(x'),
\end{equation} 
where $\theta \in (0,1)$. In the following,   we utilize the compact notation  $z = (u - u^*, x')$. \\

\noindent \emph{Quadratic bounds on the composite Lyapunov function}. Next, we derive quadratic bounds on the composite Lyapunov function. Since the function \( \phi + \psi \) is $\mu_u$-strongly convex in $u$, it holds that: 
\[ (1 - \theta)\frac{1}{\eta}\frac{\mu_{u}}{2} \|u - u^{*}\|^{2} + \frac{\theta}{\eta} \underline{\lambda} (P) \|x'\|^{2} \leq \nu (t, z; w) \]
where $u^*$ is the unique minimizer of \eqref{opt:objectiveproblem}. Moreover, since the composite function \( \phi(u) + \psi(Gu+Hw_t) \) is Lipschitz-smooth with parameter \( \ell_u + \|G\|^2 \ell_{y}  \), using the  Descent Lemma~\cite[Lemma 4.22]{Beck:2014} we get: \[ \nu (t, z; w) \leq (1 - \theta)\frac{1}{\eta} \frac{( \ell_u + \|G\|^2 \ell_{y}  )}{2}\|u - u^{*}\|^{2} + \frac{\theta}{\eta} \bar \lambda (P) \|x'\|^{2}. \] 
Setting \(c_{1}\) \(:=\min \{(1 - \theta)\frac{1}{\eta}\frac{\mu_{u}}{2},\frac{\theta}{\eta}\underline{\lambda}\left(P\right)\}\) and $ c_{2} := \max \{ (1 - \theta)\frac{1}{\eta}\frac{(\ell_u + \|G\|^2 \ell_{y} )}{2}, \frac{\theta}{\eta} \bar{\lambda} (P) \} $, obtain the quadratic bounds, 
\begin{equation}
\label{eq:quadraticBoundsLyapunov}
c_{1} \|z\|^2 \leq \nu (t, z; w) \leq c_{2} \|z\|^2.
\end{equation} 

\vspace{.4cm}

\noindent \emph{Bound for the derivatives of the Lyapunov function along the trajectories of the system}. 
Consider the derivative of the composite Lyapunov function along the trajectories of the closed loop system:
\begin{equation}
\label{eq:smallDerivativeLyapunov}
\dot \nu (t, z; w) = (1 - \theta) \frac{1}{\eta} \dot V(t, u; w) + \theta \frac{1}{\eta} \dot W(x').
\end{equation}
%
Consider the first term on the right hand side of \eqref{eq:smallDerivativeLyapunov}; the derivative can be calculated as:
\begin{align*}
    (1 - \theta) \frac{1}{\eta} \dot V(t, u; w)& =\frac{(1 - \theta)}{\eta} \frac{\partial V}{\partial u}^\top \dot u\\
    = \frac{(1 - \theta)}{\eta} \frac{\partial V}{\partial u}^\top & \eta(m(u, x' + x_{eq}) + e)\\
    = (1 - \theta) \frac{\partial V}{\partial u}^\top & (m(u, x' + x_{eq}) + e) + (1 - \theta) \frac{\partial V}{\partial u}^\top m(u,x_{eq}) - (1 - \theta) \frac{\partial V}{\partial u}^\top m(u,x_{eq})\\
    = (1 - \theta)\frac{\partial V}{\partial u}^\top & \left(m(u, x' + x_{eq}) - m(u,x_{eq})\right) + (1 - \theta) \frac{\partial V}{\partial u}^\top m(u,x_{eq}) + (1 - \theta) \frac{\partial V}{\partial u}^\top e. \\
\end{align*}
By expanding each term:
\begin{align*}
    (1 - \theta)\frac{\partial V}{\partial u}^\top &\left(m(u, x' + x_{eq}) - m(u,x_{eq})\right) =\\
    &= (1 - \theta)(\nabla \phi(u) - \nabla \phi(u^*) + G^\top \nabla \psi(Gu + Hw_t) - G^\top \nabla \psi (Gu^* + Hw_t)) \\
    & ~~ \times (-\nabla \phi(u) - G^\top \nabla \psi(Cx' + Gu + Hw_t) + \nabla \phi(u) + G^\top \nabla \psi(Gu + Hw_t))\\
    &= (1 - \theta)(\nabla \phi(u) + G^\top \nabla \psi(Gu + Hw))^\top \\
    & ~~ \times (- G^\top \nabla \psi(Cx' + Gu + Hw_t) + G^\top \nabla \psi (G u + H w_t))\\
    &\leq (1 - \theta) \|G\|\| \nabla \phi(u) + G^\top \nabla \psi(Gu + Hw_t) \| \\
    & ~~ \times \| \nabla \psi(Cx' + Gu + Hw_t) - \nabla \psi (Gu + Hw_t) \|\\
    &\leq (1 - \theta) \ell_y \|G\|\|C\|\| \nabla \phi(u) + G^\top \nabla \psi(Gu + Hw_t) \| \| x'\| , 
\end{align*}
\begin{align*}
    (1 - \theta) \frac{\partial V}{\partial u}^\top &m(u,x_{eq}) =\\
    &= (1 - \theta)(\nabla \phi(u) - \nabla \phi(u^*) + G^\top \nabla \psi(Gu + Hw_t) - G^\top \nabla \psi (Gu^* + Hw_t))^\top \\
    & ~~ \times (-\nabla \phi(u) - G^\top \nabla \psi(Gu + Hw_t))\\
    &= (1 - \theta)(\nabla \phi(u) + G^\top \nabla \psi(Gu + Hw_t))^\top \times (-\nabla \phi(u) - G^\top \nabla \psi(Gu + Hw_t))\\
    &= - (1 - \theta) \|\nabla \phi(u) + G^\top \nabla \psi(Gu + Hw_t)\|^2 , 
\end{align*}
and
\begin{align*}
    (1 - \theta) \frac{\partial V}{\partial u}^\top &e =\\
    &= (1 - \theta)(\nabla \phi(u) - \nabla \phi(u^*) + G^\top \nabla \psi(G u + H w_t) - G^\top \nabla \psi (G u^* + H w_t))^\top \\
    & ~~ \times \{(\nabla b(u) - \nabla b(u^*))^\top (\alpha - \hat \alpha_t)+ \nabla b(u^*)^\top (\alpha - \hat \alpha_t) \\
    & ~~~~ + (G^\top \nabla d(y) - G^\top \nabla d(y^*))^\top (\rho - \hat \rho_t)+ G^\top \nabla d(y^*)^\top (\rho - \hat \rho_t)\\
    & ~~~~ + (\nabla e_\phi(u^*) - \nabla e_\phi(u^*)) + (G^\top \nabla e_\psi(y) - G^\top \nabla e_\psi(y^*)) + \nabla e_\phi(u^*) + G^\top \nabla e_\psi(y^*)\}\\
    &= (1 - \theta)(\nabla \phi(u) + G^\top \nabla \psi(Gu + Hw_t))^\top \\
    & ~~ \times \{(\nabla b(u) - \nabla b(u^*))^\top (\alpha - \hat \alpha_t)+ \nabla b(u^*)^\top (\alpha - \hat \alpha_t) \\
    & ~~~~ + (G^\top \nabla d(y) - G^\top \nabla d(y^*))^\top (\rho - \hat \rho_t)+ G^\top \nabla d(y^*)^\top (\rho - \hat \rho_t)\\
    & ~~~~ + (\nabla e_\phi(u^*) - \nabla e_\phi(u^*)) + (G^\top \nabla e_\psi(y) - G^\top \nabla e_\psi(y^*)) + \nabla e_\phi(u^*) + G^\top \nabla e_\psi(y^*)\}\\
    &\leq (1 - \theta)\|\nabla \phi(u) + G^\top \nabla \psi(Gu + Hw_t)\| \\
    & ~~ \times \|(\nabla b(u) - \nabla b(u^*))^\top (\alpha - \hat \alpha_t)+ \nabla b(u^*)^\top (\alpha - \hat \alpha_t) \\
    & ~~~~ + (G^\top \nabla d(y) - G^\top \nabla d(y^*))^\top (\rho - \hat \rho_t)+ G^\top \nabla d(y^*)^\top (\rho - \hat \rho_t)\\
    & ~~~~ + (\nabla e_\phi(u^*) - \nabla e_\phi(u^*)) + (G^\top \nabla e_\psi(y) - G^\top \nabla e_\psi(y^*)) + \nabla e_\phi(u^*) + G^\top \nabla e_\psi(y^*)\|\\
    &\leq (1 - \theta) \left(\ell_u + \ell_y \|G\|^2\right)\|u - u^*\| \\
    & ~~ \times \{\|\nabla b(u) - \nabla b(u^*)\| \|\alpha - \hat \alpha_t\| + \|G\|\|\nabla d(y) - \nabla d(y^*)\| \|\rho - \hat \rho_t\| \\
    & ~~~~~~~ + \|\nabla b(u^*)\| \|\alpha - \hat \alpha_t\| + \|G\|\|\nabla d(y^*)\| \|\rho - \hat \rho_t\| \\
    & ~~~~ + \|\nabla e_\phi(u^*) - \nabla e_\phi(u^*)\| + \|G\| \|\nabla e_\psi(y) - \nabla e_\psi(y^*)\| + \|\nabla e_\phi(u^*)\| + \|G\|\|\nabla e_\psi(y^*)\|\}\\
    &= (1 - \theta) \left(\ell_u + \ell_y \|G\|^2\right)\|u - u^*\|\\
    & ~~ \times \{(\ell_u^N \|\alpha - \hat \alpha_t\| + \ell_y^M \|G\|^2 \|\rho - \hat \rho_t\| + \ell_u^e + \ell_y^e \|G\|^2)\|u - u^*\| \\
    &~~~~ + \|\nabla b(u^*)\| \|\alpha - \hat \alpha_t\| + \|G\| \| \nabla d(y^*) \| \|\rho - \hat \rho_t\|+ \|\nabla e_\phi(u^*)\| + \|G\|\|\nabla e_\psi(y^*)\|\}\\
    &= (1 - \theta) \left(\ell_u + \ell_y \| G \|^2\right)\left(\ell_u^N \|\alpha - \hat \alpha_t\| + \ell_y^M \|G\|^2 \|\rho - \hat \rho_t\| + \ell_u^e + \ell_y^e \|G\|^2\right) \|u - u^*\|^2 \\
    & ~~~ + (1 - \theta)\left(\ell_u  + \ell_y \| G \|^2 \right) \{\|\nabla b(u^*)\| \|\alpha - \hat \alpha_t\| + \|G\| \| \nabla d(y^*) \| \|\rho - \hat \rho_t\| \\
    & ~~~~~~~~~~~~~~~~~~~~~~~~~~~~~~~~~~~~~~~~~~~~~~~~~ + \|\nabla e_\phi(u^*)\| + \|G\|\|\nabla e_\psi(y^*)\|\} \|u - u^*\|.
\end{align*}

For the second term on the right hand size of \eqref{eq:smallDerivativeLyapunov}, calculate:
\begin{align*}
    \frac{\theta}{\eta} \dot W(x') &= \frac{\theta}{\eta} \frac{\partial W}{\partial x'}^\top \dot x'\\
    &= \frac{\theta}{\eta} \frac{\partial W}{\partial x'}^\top \left(n(u, x' + x_{eq}) - \eta \frac{\partial x_{eq}}{\partial u} \dot u - \frac{\partial x_{eq}}{\partial w_t} \dot w_t \right)\\
    &= \frac{\theta}{\eta} \frac{\partial W}{\partial x'}^\top \left( n(u, x' + x_{eq}) - \eta \frac{\partial x_{eq}}{\partial u} (m(u, x' + x_{eq}) + e) - \frac{\partial x_{eq}}{\partial w_t} \dot w_t \right)\\
    &= \frac{\theta}{\eta} \frac{\partial W}{\partial x'}^\top n(u, x' + x_{eq}) - \frac{\theta}{\eta} \eta \frac{\partial W}{\partial x'}^\top \frac{\partial x_{eq}}{\partial u} m(u, x' + x_{eq})  - \frac{\theta}{\eta} \eta \frac{\partial W}{\partial x'}^\top\frac{\partial x_{eq}}{\partial u} e - \frac{\theta}{\eta} \frac{\partial W}{\partial x'}^\top\frac{\partial x_{eq}}{\partial w_t} \dot w_t\\
    &= \frac{\theta}{\eta} \frac{\partial W}{\partial x'}^\top n(u, x' + x_{eq}) - \theta \frac{\partial W}{\partial x'}^\top\frac{\partial x_{eq}}{\partial u}m(u, x' + x_{eq}) - \theta \frac{\partial W}{\partial x'}^\top\frac{\partial x_{eq}}{\partial u}e - \frac{\theta}{\eta}\frac{\partial W}{\partial x'}^\top\frac{\partial x_{eq}}{\partial w_t} \dot w_t.
\end{align*}

Analyzing each term, one has that:
\begin{align*}
    \frac{\theta}{\eta} \frac{\partial W}{\partial x'}^\top n(u, x' + x_{eq}) &= - \frac{\theta}{\eta} x'^\top Q x' 
    \leq -\frac{\theta}{\eta} \underline \lambda (Q) \|x'\|^2 , 
\end{align*}
\begin{align*}
    - \theta \frac{\partial W}{\partial x'}^\top\frac{\partial x_{eq}}{\partial u}m(u, x' + x_{eq}) &= \theta (2 x'^\top P^\top)(-A^{-1}B)(\nabla \phi(u) + G^\top \nabla \psi(Cx' + Gu + Hw_t))\\
    &= \theta (2 x'^\top P^\top)(-A^{-1}B)(\nabla \phi(u) + G^\top \nabla \psi(Gu + Hw_t)\\
    & ~~~~~ + G^\top \nabla \psi(Cx' + Gu + Hw_t) - G^\top \nabla \psi(Gu + Hw_t))\\
    &\leq \theta 2 \|PA^{-1}B\|\|x'\|\{ \|\nabla \phi(u) + G^\top \nabla \psi(Gu + Hw_t)\|\\
    & ~~~~~ + \|G\|\| \nabla \psi(Cx' + Gu + Hw_t)  - \nabla \psi(Gu + Hw_t)\| \}\\
    &\leq \theta 2 \|PA^{-1}B\|\|x'\|\{ \|\nabla b(u)^\top\alpha + G^\top \nabla \psi(Gu + Hw_t)\| + \ell_y \|G\|\|C\|\| x'\| \}\\
    &= \theta 2 \|PA^{-1}B\|\|\nabla \phi(u) + G^\top \nabla \psi(Gu + Hw_t)\|\|x'\| \\
    & ~~~ + \theta 2 \|PA^{-1}B\|\ell_y \|G\|\|C\|\|x'\|^2 , 
\end{align*}
\begin{align*}
    - \theta \frac{\partial W}{\partial x'}^\top\frac{\partial x_{eq}}{\partial u}e &= \theta (2 x'^\top P^\top)(-A^{-1}B) \\
    & ~~ \times \{(\nabla b(u) - \nabla b(u^*))^\top (\alpha - \hat \alpha_t)+ \nabla b(u^*)^\top (\alpha - \hat \alpha_t) \\
    & ~~~~ + (G^\top \nabla d(y) - G^\top \nabla d(y^*))^\top (\rho - \hat \rho_t)+ G^\top \nabla d(y^*)^\top (\rho - \hat \rho_t)\\
    & ~~~~ + (\nabla e_\phi(u^*) - \nabla e_\phi(u^*)) + (G^\top \nabla e_\psi(y) - G^\top \nabla e_\psi(y^*)) + \nabla e_\phi(u^*) + G^\top \nabla e_\psi(y^*)\}\\
    &\leq \theta 2 \|PA^{-1}B\|\|x'\| \{ \ell_u^N \| \alpha - \hat \alpha_t \|\|u - u^*\| + \|\nabla b(u^*)\| \|\alpha - \hat \alpha_t\| \\
    & ~~~ + \ell_y^M \|G\|^2 \|\rho - \rho_t\|\|u - u^*\| + \|G\|\|\nabla d(y^*) \| \|\rho - \hat \rho_t\|\\
    & ~~~ + (\ell_u^e + \ell_y^e\|G\|^2)\|u-u^*\| + \|\nabla e_\phi(u)\| + \|G\| \|\nabla e_\psi(y^*)\|\} \\
    &= \theta 2 \|PA^{-1}B\|\left(\ell_u^N \|\alpha - \hat \alpha_t\| + \ell_y^M \|G\|^2 \|\rho - \hat \rho_t\| + \ell_u^e + \ell_y^e \|G\|^2\right) \|u - u^*\|\|x'\| \\
    & ~~~ + \theta 2 \|PA^{-1}B\| \{\|\nabla b(u^*)\| \|\alpha - \hat \alpha_t\| + \|G\| \| \nabla d(y^*) \| \|\rho - \hat \rho_t\| \\
    & ~~~~~~~~~~~~~~~~~~~~~~~~~~~~~~~~~~~~~~~~~~~~~~~~~ + \|\nabla e_\phi(u^*)\| + \|G\|\|\nabla e_\psi(y^*)\|\} \|x'\| , 
\end{align*}
and
\begin{align*}
    - \frac{\theta}{\eta}\frac{\partial W}{\partial x'}^\top\frac{\partial x_{eq}}{\partial w_t} \dot w_t &= \frac{1}{\eta}\theta (2 P x')^\top(-A^{-1} E) \dot w_t\\
    &\leq \frac{1}{\eta}\theta 2 \|PA^{-1}E\|\|x'\|\|\dot w_t\|.
\end{align*}
By combining the above bounds, \eqref{eq:smallDerivativeLyapunov} can be bounded above as follows:
\begin{equation}
\begin{split}
\label{eq:derivativeLyapunovUpperB}
&\dot \nu (u,w, x') \leq (1 - \theta) \{\ell_y \|G\|\|C\|\| \nabla \phi(u) + G^\top \nabla \psi(Gu + Hw_t) \| \| x'\| \\
& ~~~ - \|\nabla \phi(u) + G^\top \nabla \psi(Gu + Hw_t)\|^2\\
& ~~~ + (\ell_u + \ell_y \| G \|^2) \left(\ell_u^N \|\alpha - \hat \alpha_t\| + \ell_y^M \|G\|^2 \|\rho - \hat \rho_t\| + \ell_u^e + \ell_y^e \|G\|^2\right) \|u - u^*\|^2 \\
& ~~~ + (\ell_u + \ell_y \| G \|^2)\left(\|\nabla b(u^*)\| \|\alpha - \hat \alpha_t\| + \|G\| \| \nabla d(y^*) \| \|\rho - \hat \rho_t\|+ \|\nabla e_\phi(u^*)\| + \|G\|\|\nabla e_\psi(y^*)\|\right)\|u - u^*\|\}\\
    & - \frac{\theta}{\eta} \underline \lambda (Q) \|x'\|^2\\
    & ~~~ + \theta 2 \|PA^{-1}B\|\|\nabla \phi(u) + G^\top \nabla \psi(Gu + Hw_t)\|\|x'\| + \theta 2 \|PA^{-1}B\|\ell_y \|G\|\|C\|\|x'\|^2\\
    & ~~~ + \theta 2 \|PA^{-1}B\|\left(\ell_u^N \|\alpha - \hat \alpha_t\| + \ell_y^M \|G\|^2 \|\rho - \hat \rho_t\| + \ell_u^e + \ell_y^e \|G\|^2\right) \|x'\|\|u - u^*\| \\
    & ~~~~ + \theta 2 \|PA^{-1}B\|\left(\|\nabla b(u^*)\| \|\alpha - \hat \alpha_t\| + \|G\| \| \nabla d(y^*) \| \|\rho - \hat \rho_t\|+ \|\nabla e_\phi(u^*)\| + \|G\|\|\nabla e_\psi(y^*)\|\right)\|x'\|\\
    & ~~~~~ + \frac{1}{\eta}\theta 2 \|PA^{-1}E\|\|x'\|\|\dot w_t\|.
\end{split}
\end{equation}

\subsection{Deriving Conditions on the Controller Gain}
We will now utilize the bounds \eqref{eq:quadraticBoundW}, \eqref{eq:quadraticBoundV} to identify sufficient conditions on the controller gain to guarantee strong decrease of the Lyapunov function. 

To this aim, we will show that the bound \eqref{eq:derivativeLyapunovUpperB}  can be rewritten as
\begin{align*}
    \dot \nu &\leq - \Omega^\top \Lambda \Omega - c_0 \nu 
    +(1-\theta) \{ (\ell_u + \ell_y \| G \|^2) \left(\ell_u^N \|\alpha - \hat \alpha_t\| + \ell_y^M \|G\|^2 \|\rho - \hat \rho_t\| + \ell_u^e + \ell_y^e \|G\|^2\right) \|u - u^*\|^2\\
    & + (\ell_u + \ell_y \| G \|^2)\left(\|\nabla b(u^*)\| \|\alpha - \hat \alpha_t\| + \|G\| \| \nabla d(y^*) \| \|\rho - \hat \rho_t\|+ \|\nabla e_\phi(u^*)\| + \|G\|\|\nabla e_\psi(y^*)\|\right)\|u - u^*\|) \}\\
    & ~~~ + \theta 2 \|PA^{-1}B\|\left(\ell_u^N \|\alpha - \hat \alpha_t\| + \ell_y^M \|G\|^2 \|\rho - \hat \rho_t\| + \ell_u^e + \ell_y^e \|G\|^2\right) \|x'\|\|u - u^*\| \\
    & ~~~~ + \theta 2 \|PA^{-1}B\|\left(\|\nabla b(u^*)\| \|\alpha - \hat \alpha_t\| + \|G\| \| \nabla d(y^*) \| \|\rho - \hat \rho_t\|+ \|\nabla e_\phi(u^*)\| + \|G\|\|\nabla e_\psi(y^*)\|\right)\|x'\|\\
    & ~~~~~ + \frac{1}{\eta}\theta 2 \|PA^{-1}E\|\|x'\|\|\dot w_t\|,
\end{align*}
where $\Lambda$ is a $2 \times 2$ positive definite matrix and 
\begin{equation}\label{eq:VectorEqn}
\begin{split}
\Omega(x^\prime,u) = \begin{bmatrix} 
\|\nabla \phi(u) + G^\top \nabla \psi(Gu + H w_t)\| \\
    \|x'\|\end{bmatrix}.    
\end{split}
\end{equation}
Notice that, in this case, $\Omega(x^\prime,u)^\top \Lambda \Omega(x^\prime,u) = 0$ if and only if $(x^\prime,u)=(0,u^*_t)$ (namely, at the optimizer) and $\Omega(x^\prime,u)^\top \Lambda \Omega(x^\prime,u) > 0$
otherwise.

To this aim, let $s \in (0, 1)$, and re-organize~\eqref{eq:derivativeLyapunovUpperB} as: 
\begin{equation}
\begin{split}
\label{eq:derivativeLyapunovUpperB2}
&\dot \nu (u,w, x') \leq (1 - \theta) \{\ell_y \|G\|\|C\|\| \nabla \phi(u) + G^\top \nabla \psi(Gu + Hw_t)\| \| x'\| \\
& ~~~ - (1 - s) \|\nabla \phi(u) + G^\top \nabla \psi(Gu + Hw_t)\|^2 - s \|\nabla \phi(u) + G^\top \nabla \psi(Gu + Hw_t)\|^2\\
& ~~~ + (\ell_u  + \ell_y \| G \|^2) \left(\ell_u^N \|\alpha - \hat \alpha_t\| + \ell_y^M \|G\|^2 \|\rho - \hat \rho_t\| + \ell_u^e + \ell_y^e \|G\|^2\right) \|u - u^*\|^2 \\
& ~~~ + (\ell_u + \ell_y \| G \|^2)\left(\|\nabla b(u^*)\| \|\alpha - \hat \alpha_t\| + \|G\| \| \nabla d(y^*) \| \|\rho - \hat \rho_t\|+ \|\nabla e_\phi(u^*)\| + \|G\|\|\nabla e_\psi(y^*)\|\right)\|u - u^*\|\}\\
    & - (1 - s)\frac{\theta}{\eta} \underline \lambda (Q) \|x'\|^2 - s \frac{\theta}{\eta} \underline \lambda (Q) \|x'\|^2\\
    & ~~~ + \theta 2 \|PA^{-1}B\|\|\nabla \phi(u) + G^\top \nabla \psi(Gu + Hw_t)\|\|x'\| + \theta 2 \|PA^{-1}B\|\ell_y \|G\|\|C\|\|x'\|^2\\
    & ~~~ + \theta 2 \|PA^{-1}B\|\left(\ell_u^N \|\alpha - \hat \alpha_t\| + \ell_y^M \|G\|^2 \|\rho - \hat \rho_t\| + \ell_u^e + \ell_y^e \|G\|^2\right) \|x'\|\|u - u^*\| \\
    & ~~~~ + \theta 2 \|PA^{-1}B\|\left(\|\nabla b(u^*)\| \|\alpha - \hat \alpha_t\| + \|G\| \| \nabla d(y^*) \| \|\rho - \hat \rho_t\|+ \|\nabla e_\phi(u^*)\| + \|G\|\|\nabla e_\psi(y^*)\|\right)\|x'\|\\
    & ~~~~~ + \frac{1}{\eta}\theta 2 \|PA^{-1}E\|\|x'\|\|\dot w_t\|.
\end{split}
\end{equation}
In particular, consider the terms \( s \|\nabla \phi(u) + G^\top \nabla \psi(Gu + Hw_t)\|^2\) and \( s \frac{\theta}{\eta} \underline \lambda (Q) \|x'\|^2\) in \eqref{eq:derivativeLyapunovUpperB2}. For the first term, utilize the PL inequality to obtain, 
\begin{align*}
    - s(1 - \theta) &\|\nabla \phi(u) + G^\top \nabla \psi(Gu + Hw_t)\|^2 \leq\\
    & ~~~ - s 2(1 - \theta) \mu_u \| \phi(u)  - \phi(u^*) + \psi(Gu + Hw_t) - \psi(Gu^* + Hw_t)\|\\
    &= - s 2(1 - \theta) \mu_u (\phi(u)  - \phi(u^*) + \psi(Gu + Hw_t) - \psi(Gu^* + Hw_t)\\
    &= - s 2(1 - \theta) \mu_u V.
\end{align*}
For the second term,  calculate \[ \dot W(x') = -x'^\top Q x' \leq - \frac{\underline{\lambda}(Q)}{\bar \lambda (P)} \|x'\|^2 = - \frac{\underline{\lambda}(Q)}{\bar \lambda (P)} W(x'). \]
Then the second term simplifies to, 
\begin{align*}
    -s \frac{\theta}{\eta}\underline{\lambda}(Q)\|x'\|^2 &\leq -s \frac{\theta}{ \eta}\frac{\underline{\lambda}(Q)}{\bar \lambda (P)}W(x').
\end{align*}
Together, these give, 
\begin{equation}\label{eq:recoverNu}
    \begin{split}
        - s(1 - \theta) &\|\nabla \phi(u) + G^\top \nabla \psi(Gu + Hw_t)\|^2 -s \frac{\theta}{ \eta}\underline{\lambda}(Q)\|x'\|^2\\
        & \leq - s (1 - \theta)2 \mu_u V - s \frac{\theta}{ \eta}\frac{\underline{\lambda}(T)}{\bar \lambda (P)}W\\
        &\leq \min \left\{ s 2\mu_u\eta ,  s \frac{\underline{\lambda}(Q)}{\bar \lambda (P)}\right\} \nu := c_0 \nu.
    \end{split}
\end{equation}


To obtain the quadratic form $\Omega(x^\prime,u)^\top \Lambda \Omega(x^\prime,u)$, define the following coefficients: 
\begin{align*}
    \alpha_1 &:= (1 - s),\\
    \alpha_2 &:= (1 - s)\underline{\lambda}(Q),\\
    \beta_1 &:= \ell_y \|G\|\|C\|,\\
    \beta_2 &:= 2\|PA^{-1}B\|.\\
\end{align*}
Then, using  terms in~\eqref{eq:derivativeLyapunovUpperB2}, one can build $\Lambda$ as
\begin{align*}
    \Lambda &=\left(\begin{array}{ll}(1 - \theta)\alpha_1 & -\frac{1}{2}[(1-\theta)\beta_1 + \theta \beta_2] \\ -\frac{1}{2}[(1-\theta)\beta_1 + \theta \beta_2] & \theta[\frac{\alpha_2}{\eta} - \beta_1\beta_2]
    \end{array}\right) \\
\end{align*}
Matrix $\Lambda$ is positive definite if and only if its principal minors are positive. This yields the following condition on $\eta$:
\begin{align*}
    [(1-\theta)\alpha_1][\theta \left(\frac{\alpha_2}{\eta}- \beta_1  \beta_2\right)] &> \frac{1}{4}[(1-\theta)\beta_1 + \theta \beta_2]^2\\
    \alpha_1[\frac{\alpha_2}{\eta} - \beta_1  \beta_2] &> \frac{1}{4\theta (1-\theta)}[(1-\theta)\beta_1 + \theta \beta_2]^2\\
    \alpha_1\frac{\alpha_2}{\eta} &> \frac{1}{4\theta (1-\theta)}[(1-\theta)\beta_1 + \theta \beta_2]^2 + \alpha_1 \beta_1 \beta_2,\\
\end{align*}
and, thus, 
\begin{equation}\label{eq:etaBound}
    \eta < \frac{ \alpha_1\alpha_2}{\frac{1}{4\theta (1-\theta)}[(1-\theta)\beta_1 + \theta \beta_2]^2 + \alpha_1 \beta_1  \beta_2} .  
\end{equation}
The right hand side of \eqref{eq:etaBound} is a concave function of $\theta$ and, by maximizing with respect to $\theta$, we have that the maximum is obtained at 
$$\theta^* = \frac{\beta_1}{\beta_1 + \beta_2},$$ which gives $$\eta^* = \frac{\alpha_1 \alpha_2}{(\beta_1 \beta_2)(1 + \alpha_1)}.$$

\subsection{Analysis of the Learning Error and Derivation of the Main Result }

Finally, we will show that the time-derivative of the Lyapunov function can be bounded as
\begin{equation*}
    \dot \nu \leq - \Omega(x^\prime,u)^\top \Lambda \Omega(x^\prime,u) - c_0 \nu + c_3 \gamma(t) \nu + \left(c_4 \delta(t) + c_5 \|\dot w_t\| \right)\sqrt{\nu}. 
\end{equation*}

To this aim, recall that $z := (u - u^*,x')$. We will use the inequality $\|u - u^*\| + \|x'\| \leq \sqrt{2}\|z\|$.

We begin by using \eqref{eq:quadraticBoundV}-\eqref{eq:quadraticBoundW}, 
to obtain 
\( \|u - u^*\|^2 \leq \frac{2}{\mu_u}V\)
and 
$\|x'\|^2 \leq \frac{W}{\underline{\lambda}(P)}$. 
Further, by using 
$\nu := (1 - \theta) \frac{1}{\eta}V + \frac{\theta}{\eta}W,$ we have that \( V \leq \frac{\eta}{(1 - \theta)}\nu\) and \( \sqrt{V} \leq \frac{\sqrt{\eta}}{\sqrt{(1 - \theta)}} \sqrt{\nu}.\) Using these facts, we can calculate the following bounds:

\noindent (a)
\begin{align}
    (1 - \theta)(\ell_u  + \ell_y \| G \|^2)& \left(\ell_u^N \|\alpha - \hat \alpha_t\| + \ell_y^M \|G\|^2 \|\rho - \hat \rho_t\| + \ell_u^e + \ell_y^e \|G\|^2\right) \|u - u^*\|^2  \nonumber \\
    & \hspace{-2.5cm}\leq \frac{2}{\mu_u} (1 - \theta)(\ell_u  + \ell_y \| G \|^2) \left(\ell_u^N \|\alpha - \hat \alpha_t\| + \ell_y^M \|G\|^2 \|\rho - \hat \rho_t\| + \ell_u^e + \ell_y^e \|G\|^2\right) V(u) \label{eq:a}
\end{align}
\noindent (b) with $\kappa := (1 - \theta)(\ell_u + \ell_y \| G \|^2)(\|\nabla b(u^*)\| \|\alpha - \hat \alpha_t\| $ $+ \|G\| \| \nabla d(y^*) \| \|\rho - \hat \rho_t\|+ \|\nabla e_\phi(u^*)\| + \|G\|\|\nabla e_\psi(y^*)\|)$ for brevity, 
\begin{align}
    \kappa \|u - u^*\| \leq \kappa \sqrt{\frac{2}{\mu_u}}\sqrt{V(u)}. \label{eq:b}
\end{align}
\noindent (c)
\begin{align}
    2 \theta \|PA^{-1}B\|& \left(\ell_u^N \|\alpha - \hat \alpha_t\| + \ell_y^M \|G\|^2 \|\rho - \hat \rho_t\| + \ell_u^e + \ell_y^e \|G\|^2\right)\|x'\|\|u - u^*\| \nonumber \\ 
    & \leq 4 \theta  \|PA^{-1}B\| \left(\ell_u^N \|\alpha - \hat \alpha_t\| + \ell_y^M \|G\|^2 \|\rho - \hat \rho_t\| + \ell_u^e + \ell_y^e \|G\|^2\right)\|z\|^2 \nonumber \\
    &\leq 4 \frac{\theta}{c_1} \|PA^{-1}B\| \left(\ell_u^N \|\alpha - \hat \alpha_t\| + \ell_y^M \|G\|^2 \|\rho - \hat \rho_t\| + \ell_u^e + \ell_y^e \|G\|^2\right) \nu \nonumber \\
    &\leq \frac{4}{c_1} \|PA^{-1}B\| \left(\ell_u^N \|\alpha - \hat \alpha_t\| + \ell_y^M \|G\|^2 \|\rho - \hat \rho_t\| + \ell_u^e + \ell_y^e \|G\|^2\right) \nu  \label{eq:c}
\end{align}
since $\theta \in (0,1)$. 

\noindent (d) with $\kappa^\prime := 2 \theta  \|PA^{-1}B\|(\|\nabla b(u^*)\| \|\alpha - \hat \alpha_t\|$ $+ \|G\| \| \nabla d(y^*) \| \|\rho - \hat \rho_t\|+ \|\nabla e_\phi(u^*)\| + \|G\|\|\nabla e_\psi(y^*)\|)$ for brevity, 
\begin{align}
    \kappa^\prime \|x'\| \leq \kappa^\prime \frac{\sqrt{W}}{\sqrt{\underline{\lambda}(P)}}.\label{eq:d}
\end{align}
\noindent And, (e)
\begin{align}
    \frac{\theta}{\eta}2\|P^\top A^{-1}E\| \|\dot w_t\| \|x'\| \leq \frac{\theta}{\eta}2\|P^\top A^{-1}E\| \|\dot w_t\| \frac{\sqrt{W}}{\sqrt{\underline{\lambda}(P)}}. \label{eq:e}
\end{align} 
Define 
\begin{align}
    \gamma(t) & := \ell_u^N \|\alpha - \hat \alpha_t\| + \ell_y^M \|G\|^2 \|\rho - \hat \rho_t\| + \ell_u^e + \ell_y^e \|G\|^2, \nonumber\\
    \delta(t) &:= \|\nabla b(u^*)\| \|\alpha - \hat \alpha_t\| + \|G\| \| \nabla d(y^*) \| \|\rho - \hat \rho_t\|+ \|\nabla e_\phi(u^*)\| + \|G\|\|\nabla e_\psi(y^*)\|.
\end{align}
For~\eqref{eq:a}, obtain a bound with $\nu$ as: 
\begin{align*}
    (1 - \theta)(\ell_u  + \ell_y \| G \|^2) \gamma(t) \frac{2}{\mu_u}V &\leq (\ell_u  + \ell_y \| G \|^2) \gamma(t) \frac{2\eta}{\mu_u}\nu.
\end{align*}
From the second inequality~\eqref{eq:b}, similarly obtain:
\begin{align*}
    (1 - \theta)(\ell_u + \ell_y \| G \|^2)\delta(t)\sqrt{\frac{2}{\mu_u}}\sqrt{V} &\leq \sqrt{(1 - \theta)}(\ell_u + \ell_y \| G \|^2)\delta(t)\sqrt{\frac{2}{\mu_u}}\sqrt{\eta} \sqrt{\nu} \\
    &\leq (\ell_u + \ell_y \| G \|^2)\delta(t) \sqrt{\frac{2}{\mu_u}}\sqrt{\eta}\sqrt{\nu}.
\end{align*}
The third term remains as-is because it is already bounded above by $\nu$. 
For the fourth term, obtain: 
\begin{align*}
    2 \theta \|PA^{-1}B\|\delta(t)\frac{\sqrt{W}}{\sqrt{\underline{\lambda}(P)}} &\leq 2 \sqrt{\theta } \|PA^{-1}B\|\delta(t)\frac{\sqrt{\nu}}{\sqrt{\underline{\lambda}(P)}} \sqrt{\eta}\\
    &\leq 2 \|PA^{-1}B\|\delta(t) \sqrt{\nu}\frac{\sqrt{\eta}}{\sqrt{\underline{\lambda}(P)}}.
\end{align*}
Similarly, we bound the term due to the time-varying disturbance $w_t$ as, 
\begin{align*}
    \frac{\theta}{\eta}2\|P^\top A^{-1}E\| \|\dot w_t\| \frac{\sqrt{W}}{\sqrt{\underline{\lambda}(P)}} &\leq \frac{\theta}{\eta}2\|P^\top A^{-1}E\| \|\dot w_t\| \frac{\sqrt{\eta}}{\sqrt{\underline{\lambda}(P)}} \sqrt{\nu }\\
    &\leq \frac{2 \|P^\top A^{-1} E\|}{\sqrt{\eta} \sqrt{\underline{\lambda}(P)}}\| \dot w_t\| \sqrt{\nu}.
\end{align*}
Excluding $\gamma(t)$ and $\delta(t)$, group the terms above in front of $\eta$ and $\sqrt{\eta}$ and define: \[ c_3 := \max \left\{ (\ell_u  + \ell_y \| G \|^2) \frac{2\eta}{\mu_u}, \|P A^{-1}B\| \frac{4}{c_1}\right\}\] \[\text{ and } c_4 := \sqrt{\eta} \max \left\{ (\ell_u + \ell_y \| G \|^2) \sqrt{\frac{2}{\mu_u}}, \|PA^{-1}B\| \frac{2}{\sqrt{\underline{\lambda}(P)}}\right\}.\]
For the term dealing with $\| \dot w_t \|$, similarly define:
\[ c_5 := \frac{2 \|P^\top A^{-1} E\|}{\sqrt{\eta} \sqrt{\underline{\lambda}(P)}}. \]
In summary, 
\begin{equation}
    \begin{split}
        \gamma(t) &:= \ell_u^N \|\alpha - \hat \alpha_t\| + \ell_y^M \|G\|^2 \|\rho - \hat \rho_t\| + \ell_u^e + \ell_y^e \|G\|^2, \\
        \delta(t) &:= \|\nabla b(u^*)\| \|\alpha - \hat \alpha_t\| + \|G\| \| \nabla d(y^*) \| \|\rho - \hat \rho_t\|+ \|\nabla e_\phi(u^*)\| + \|G\|\|\nabla e_\psi(y^*)\|, \\
        c_0 &:= s \min \left\{ 2 \mu_u \eta, \frac{\underline{\lambda}(Q)}{\bar 2 \lambda (P)} \right\}, \\
        c_1 &:= \min \left\{(1 - \theta)\frac{1}{\eta}\frac{\mu_{u}}{2},\frac{\theta}{\eta}\underline{\lambda}\left(P\right)\right\}, \\
        c_2 &:= \max \left\{ (1 - \theta)\frac{1}{\eta}\frac{(\ell_u + \|G\|^2 \ell_{y} )}{2}, \frac{\theta}{\eta} \bar{\lambda} (P) \right\}, \\
        c_3 &:= \max \left\{ (\ell_u  + \ell_y \| G \|^2) \frac{2\eta}{\mu_u}, \|P A^{-1}B\| \frac{4}{c_1}\right\}, \\
        c_4 &:= \sqrt{\eta} \max \left\{ (\ell_u + \ell_y \| G \|^2) \sqrt{\frac{2}{\mu_u}}, \|PA^{-1}B\| \frac{2}{\sqrt{\underline{\lambda}(P)}}\right\}, \\
        c_5 &:= \frac{2 \|P^\top A^{-1} E\|}{\sqrt{\eta} \sqrt{\underline{\lambda}(P)}}.
    \end{split}
\end{equation}
Altogether, our analysis bounds $\dot \nu$ as, 
\begin{equation}
\label{eq:bound_dotnu}
    \dot \nu \leq - \Omega(x^\prime,u)^\top \Lambda \Omega(x^\prime,u) - c_0 \nu + c_3 \gamma(t) \nu + \left(c_4 \delta(t) + c_5 \|\dot w_t\| \right)\sqrt{\nu}. 
\end{equation}
We have derived sufficient conditions on $\eta$ so that $\Lambda$ is positive definite; consequently, $\Omega(x^\prime,u)^\top \Lambda \Omega(x^\prime,u) > 0$ must hold and~\eqref{eq:bound_dotnu} can be further bounded as: 
\begin{equation}\label{eq:lastStep}
    \dot \nu \leq - \left(c_0  - c_3 \gamma(t) \right) \nu + \left(c_4 \delta(t) + c_5 \|\dot w_t\|\right) \sqrt{\nu}. 
\end{equation}
The following analysis  for \eqref{eq:lastStep} is inspired by \cite[Section 9.3]{Khalil:1173048}. 
First, apply a change of variables $\nu' := \sqrt{\nu}$. Then, 
\[ \dot \nu' \leq -\frac{1}{2}\left( c_0 - c_3 \gamma (t) \right) \nu' + \frac{1}{2}\left( c_4 \delta(t) + c_5 \|\dot w_t\| \right). \]
Let $\Phi(t,t_0) := e^{-\frac{1}{2}c_0(t - t_0) + \frac{c_3}{2}\int_{t_0}^t \gamma(\tau) d\tau },$ for a given $t_0 \geq 0$.
By the Comparison Lemma \cite[Section 3.4, Section 9.3]{Khalil:1173048}, 
\[ \dot \nu' \leq \Phi(t, t_0) \nu'(t_0)  + \frac{c_4}{2} \int_{t_0}^t \Phi(t, \tau) \delta(\tau) d\tau + \frac{c_5}{2}\int_{t_0}^t \Phi(t, \tau)  \|\dot w_{\tau}\|d\tau. \]
By \eqref{eq:quadraticBoundsLyapunov}, 
\begin{equation}
    \|z(t)\| \leq \sqrt{\frac{c_2}{c_1}}\Phi(t,t_0) \|z(t_0)\| + \frac{c_4}{2 \sqrt{c_1}} \int_{t_0}^t \Phi(t, \tau)\delta(\tau) d\tau + \frac{c_5}{2 \sqrt{c_1}} \int_{t_0}^t \Phi(t, \tau)\|\dot w_{\tau}\| d\tau.
\end{equation}

Note that the errors $\|\alpha - \hat \alpha_t\|$ and $\|\rho - \hat \rho_t\|$ within $\gamma(t)$ and $\delta(t)$ are piece-wise constant.
Thus, for any $K, K' \in \naturalset$,  we bound $\gamma(t)$ as: 
\begin{align}
    \sum_{i=0}^K & \ell_u^N \|\alpha - \hat \alpha_{t_{i-1}}\|(t_i - t_{i-1})  + \sum_{i = 0}^{K'} \ell_y^M \|G\|^2 \|\rho - \hat \rho_{t_{i-1}}\|   (t_i - t_{i-1}) \nonumber \\
    &\leq \underbrace{\left(\ell_u^N \sup_{t_0 \leq t \leq \tau} \left\{ \|\alpha - \hat \alpha_\tau\|\right \} + \ell_y^M \|G\|^2 \sup_{t_0 \leq \tau \leq t} \left\{ \|\rho - \hat \rho_\tau\| \right\} + \ell_u^e + \ell_y^e \|G\|^2 \right)}_{:= \epsilon'}(t - t_0).  \label{eq:bounderror}
\end{align}
In \eqref{eq:lastStep}, we must ensure that \( - c_0 \nu + c_3 \gamma(t) \nu < 0\). Using~\eqref{eq:bounderror}, we impose the condition  \[\epsilon' < \frac{c_0}{c_3}. \]
Finally, defining \( a':= \frac{1}{2}\left( c_0 - \epsilon' c_3 \right) > 0\), we obtain, 
\begin{equation}
\label{eq:finalresultproof}
    \|z(t)\| \leq \sqrt{\frac{c_2}{c_1}} \|z(t_0)\| e^{-a'(t - t_0)} + \frac{c_4}{2 \sqrt{c_1}} \int_{t_0}^t e^{-a'(t - \tau)} \delta(\tau) d \tau  + \frac{c_5}{2 \sqrt{c_1}} \int_{t_0}^t e^{-a'(t - \tau)} \| \dot w_\tau \| d \tau. 
\end{equation}
Theorem~\ref{thm:truncation} follows by noting that 
$\Xi(t) = \delta(t)$.

To show the result of  Theorem~\ref{thm:finiteBasisFunctions}, the same steps can be used upon setting $\nabla e_\phi(u) = 0$, $\nabla e_\psi(y) = 0$, $\nabla e_\phi(u^*) = 0$, and $\nabla e_\psi(y^*) =0 $ in~\eqref{eq:nominalperturbed}, since there is no truncation error in Theorem~\ref{thm:finiteBasisFunctions}. This lead to the re-definition of
$\gamma(t)$ and $\delta(t)$ as $\gamma(t) := \ell_u \|\alpha - \hat \alpha_t\| + \ell_y \|G\|^2 \|\rho - \hat \rho_t\|$ and 
$\delta(t) := \|\nabla b(u^*)\| \|\alpha - \hat \alpha_t\| + \|G\| \| \nabla d(y^*) \| \|\rho - \hat \rho_t\|$, respectively. 

\section{Numerical Values for Simulation}
\label{ap:simulations}
Here, we provide the exact matrices used to generate Figure~2.

\begin{align*}
    A &= \begin{bmatrix}
    -2.7527 & -0.6944 & -2.8952 & -0.7989\\
    1.2008 & -4.3397 & -1.7097 & -0.6025\\
    -0.2198 & -1.0665 & -5.1494 & 0.3043\\
    -2.8886 & 1.922 & 2.7361 & -3.8897
    \end{bmatrix},
    B = C = D = E = \begin{bmatrix}
    1 & 0 & 0 & 0\\
    0 & 1 & 0 & 0\\
    0 & 0 & 1 & 0\\
    0 & 0 & 0 & 1
    \end{bmatrix}, \\
    \Upsilon &= \begin{bmatrix}
    6.095 & 0.6234 & 0.1468 & -0.9387\\
    0.6234 & 6.4595 &  -1.0145 &  1.0203\\
    0.1468 & -1.0145 & 7.0719 & 0.7042\\
    -0.9387 & 1.0203 & 0.7042 & 4.5038
    \end{bmatrix}, 
    \upsilon = \begin{bmatrix}
    0.0201\\
    1.4908\\
    1.2373\\
    1.8092
    \end{bmatrix},
    r = -0.1504,\\
    Q &= \begin{bmatrix}
    3.994 & -1.1602 & -0.1978 & -0.9408\\
    -1.1602 & 4.0145 & -0.3114 & -0.8189\\
    -0.1978 & -0.3114 & 5.9914 & -1.8039\\
    -0.9408 & -0.8189 & -1.8039 & 5.3419
    \end{bmatrix}, 
    P = \begin{bmatrix}
    1.3220    & 0.3400   & -0.3819   & -0.9667\\
    0.3400 &    0.7413 &  -0.3188  & -0.4261\\
   -0.3819  & -0.3188  &  0.6745  &  0.1771\\
   -0.9667  & -0.4261  &  0.1771  &  1.3186\\
    \end{bmatrix}.
\end{align*}
The matrices $Q,P$ are used in the Lyapunov function for the boundary layer system as shown in Appendix \ref{S:proofs}. 

\end{document}